\documentclass[12pt]{amsart}

\usepackage{amssymb}
\usepackage[cmtip,matrix,arrow]{xy}
\newtheorem{theorem}{Theorem}[section]
\newtheorem{corollary}[theorem]{Corollary}
\newtheorem{proposition}[theorem]{Proposition}

\newtheorem{lemma}[theorem]{Lemma}
\theoremstyle{definition}    
\newtheorem{definition}[theorem]{Definition}
\theoremstyle{remark}

\newtheorem{remark}[theorem]{Remark}
\newtheorem{remarks}[theorem]{Remarks}
\newtheorem{example}[theorem]{Example}
\newcommand{\beq}{\begin{eqnarray*}}
\newcommand{\eeq}{\end{eqnarray*}}
\renewcommand{\H}{{\mathcal H}}
\newcommand{\D}{\mathfrak{D}}

\newcommand{\W}{\mathcal{W}}

\newcommand{\La}{\Lambda}

\newcommand{\ca}{\mathcal}

\newcommand{\R}{\mathbb{R}}

\newcommand{\SU}{\on{SU}}

\newcommand{\Z}{\mathbb{Z}}

\DeclareMathOperator{\Cl}{Cl}
\newcommand\lie[1]{\mathfrak{#1}}

\newcommand{\g}{\lie{g}}

\newcommand{\on}{\operatorname}

\newcommand{\Ad}{ \on{Ad} }
\newcommand{\ad}{\on{ad}}
\newcommand{\Pf}{\on{Pf}}  

\renewcommand{\ker}{ \on{ker}}

\newcommand{\Mult}{  \on{Mult}}

\newcommand\qu{/\kern-.7ex/} 

\newcommand{\rr}{\mathfrak{r}}

\newcommand{\hra}{\hookrightarrow}

\renewcommand{\d}{{\mbox{d}}}

\newcommand\Phinv{\Phi^{-1}}

\newcommand{\End}{\on{End}}
\newcommand{\Der}{\on{Der}}

\newcommand\sig{\sigma}
\newcommand\eps{\epsilon}
\newcommand\Om{\Omega}
\newcommand\om{\omega}

\newcommand{\f}{\frac}
\newcommand{\Gr}{\on{Gr}}
\newcommand{\olt}{\overline{\theta}}
\newcommand{\p}{\partial}
\renewcommand{\l}{\langle}
\renewcommand{\r}{\rangle}
\newcommand\hh{{\f{1}{2}}}
\newcommand{\ti}{\tilde}

\newcommand{\wh}{\widehat}

\begin{document}

\title{The non-commutative Weil Algebra}
\author[A. Alekseev]{A. Alekseev$^{*}$}
\address{Institute for Theoretical Physics \\ Uppsala University \\
Box 803 \\ \mbox{S-75108} Uppsala \\ Sweden}
\email{alekseev@teorfys.uu.se}
\thanks{
$^{*}$ Partially supported by NFR under contract FU 11672-302}

\author[E. Meinrenken]{E. Meinrenken$^{\dagger}$}
\address{University of Toronto, Department of Mathematics,
100 St George Street, Toronto, Ontario M5R3G3, Canada }
\email{mein@math.toronto.edu}
\thanks{$^{\dagger}$ Partially supported by NSERC grant \# 72011899
and by a Connaught grant.}

\begin{abstract}
For any compact Lie group $G$, together with an invariant inner
product on its Lie algebra $\g$, we define the non-commutative Weil
algebra $\W_G$ as a tensor product of the universal enveloping algebra
$U(\g)$ and the Clifford algebra $\on{Cl}(\g)$. Just as the usual Weil
algebra $W_G=S\g^*\otimes \wedge\g^*$, $\W_G$ carries the structure of
an acyclic, locally free $G$-differential algebra and can be used to
define equivariant cohomology $\ca{H}_G(B)$ for any $G$-differential
algebra $B$. The main result of this paper is the construction of an
isomorphism $\ca{Q}:\,W_G\to \W_G$ of the two Weil algebras as
$G$-differential spaces.  Furthermore, we prove that the corresponding
vector space isomorphism from the usual equivariant cohomology
$H_G(B)$ to the equivariant cohomology $\H_G(B)$ is in fact a ring
isomorphism.  This generalizes the Duflo isomorphism $(S\g)^G\cong
U(\g)^G$ between the ring of invariant polynomials and the ring of
Casimir elements. We extend our considerations to Weil algebras and
equivariant cohomology with generalized coefficients, where the
algebra $U(\g)$ is replaced by the convolution algebra $\ca{E}'(G)$ of
distributions on $G$.
\end{abstract}

\maketitle

\section{Introduction}
Let $G$ be a connected Lie group with Lie algebra $\g$.  The Duflo map
is a vector space isomorphism $\on{Duf}:\,S(\g)\to U(\g)$ between the
symmetric algebra and the universal enveloping algebra which, as
proved by Duflo \cite{du:op}, restricts to a {\em ring} isomorphism 
from the algebra of invariant
polynomials $S(\g)^G$ onto the center $U(\g)^G$ of the universal
enveloping algebra.  For semi-simple $\g$ the Duflo map coincides
with the Harish-Chandra isomorphism. 
The Duflo map extends to a map of compactly
supported distributions, $\on{Duf}:\,\ca{E}'(\g)\to \ca{E}'(G)$, which
is a ring homomorphism for the $G$-invariant parts.  
For generalizations of Duflo's theorem 
see the papers of Kashiwara-Vergne \cite{ka:ca} and 
Kontsevich \cite{ko:de}.  

In this paper we will obtain analogues of the Duflo
map and of Duflo's theorem in the context of equivariant
cohomology of $G$-manifolds. Our result will involve a non-commutative
version of the de Rham model of equivariant cohomology. Throughout we
will assume that the group $G$ is compact. Let 
$$W_G=S\g^*\otimes \wedge\g^*$$
the Weil algebra. 
Given a basis $e_a$ of $\g$ let $f_{ab}^c$ the
structure constants of $\g$, and denote by $v^a,\,y^a$ the generators
of $S\g^*$ and $\wedge\g^*$ corresponding to the dual basis.  The
coadjoint action of $G$ on $\g^*$ induces an action on $W_G$, with
generators $L_a$. Let derivations $\iota_a$ be given on generators by
$\iota_a y^b=\delta_a^b$ and $\iota_a v^b=0$. In \cite{ca:no},
H. Cartan shows that the derivation
\begin{equation}\label{WDiff} 
\d y^a=v^a-\hh f^a_{jk}  y^j y^k,\ \
\d v^a=-f^a_{jk} y^j v^k
\end{equation}
gives $W_G$ the structure of a $G$-differential algebra. In
particular, $\d$ is a differential and $[\iota_a,\d]=L_a$.  

Suppose $\g$ comes equipped with an invariant inner product, used 
to identify $\g\cong\g^*$. Let $\on{Cl}(\g)$ be the corresponding Clifford 
algebra and let the {\em non-commutative Weil algebra} be its tensor  
product with the universal enveloping algebra:
$$ \W_G=U(\g)\otimes \Cl(\g).$$
Take the basis $e_a$ to be orthonormal, and let $u_a,x_a$ be the
corresponding generators of $U(\g)$ and $\Cl(\g)$. Again $\W_G$ 
carries a $G$-action induced by the adjoint action; let $L_a$ 
be its generators and let $\iota_a$ be the derivation extension of 
$\iota_a x_b=\delta_{ab}$ and $\iota_a u_b=0$. We show that there 
exists a derivation $\d$ on $\W_G$ which on generators is given 
by formulas analogous to \eqref{WDiff}: 
\begin{equation}\label{WDiff1}
\d x_a=u_a-\hh f_{ajk}  x_j x_k,\ \
\d u_a=-f_{ajk} x_j  u_k.
\end{equation}
Moreover $\W_G$ is still a $G$-differential algebra, that is,
$\d$ squares to zero and Cartan's formula continues to hold. A new feature 
is that the derivations $\iota_a,L_a,\d$ can be written as 
commutators. In  particular, $\d=\ad(\D)$ where 
$$ \D=u_a x_a -\f{1}{6} f_{abc} x_a x_b x_c$$
squares to the quadratic Casimir element 
$$ \D^2=\hh u_a u_a -\f{1}{48} f_{abc}f_{abc}.$$
The latter is naturally interpreted as a Laplace operator on 
$G$ and $\D$ as a Dirac operator. 

The main result of this paper is the construction of a vector space 
isomorphism $\ca{Q}:\, W_G\to \W_G$ (called the quantization map) 
which intertwines the derivations $\iota_a,L_a,\d$. Put differently, 
$\ca{Q}$ induces a second, non-commutative algebra structure on 
the Weil algebra for which $\iota_a,L_a,\d$ continue to be 
derivations. On the symmetric algebra $S\g\otimes 1$, $\ca{Q}$ 
restricts to the Duflo map 
while on the exterior algebra $1\otimes \wedge\g$, it restricts 
to the inverse of the symbol map $\sig:\Cl(\g)\to \wedge\g$.  
However, $\ca{Q}$ is not the direct product $\on{Duf}\otimes 
\sig^{-1}$ but has the more complicated form
$$\ca{Q}=(\on{Duf}\otimes \sig^{-1})\circ \exp(-\hh T_{ab}\iota_a\iota_b)$$
where $T_{ab}$ is a certain solution of the classical dynamical Yang-Baxter 
equation \cite{et:ge}. 

We also consider Weil algebras with generalized coefficients,
$\widehat{W}_G=\ca{E}'(\g^*)\otimes \wedge\g^*$ and
$\widehat{\W}_G=\ca{E}'(G) \otimes \on{Cl}(\g)$. Here $S\g^*$ is
identified with the subalgebra of $\ca{E}'(\g^*)$ consisting of
distributions with support at the origin and $U(\g)$ with the
subalgebra of $\ca{E}'(G)$ consisting of distributions with support at
the group unit.  The quantization map extends to a map
$\ca{Q}:\widehat{W}_G\to \widehat{\W}_G$ which still intertwines the
derivations $\iota_a,L_a,\d$ but is no longer an isomorphism.

Given any $G$-differential algebra $B$, for example the de Rham
complex $\Om^\star(M)$ of a $G$-manifold $M$, the equivariant 
cohomology $H_G(B)$ is defined as 
the cohomology algebra of the
basic subcomplex of $W_G\otimes B$. It is a module over 
the ring $(S\g^*)^G$ of invariant polynomials. Replacing $W_G$ with 
$\wh{W}_G,\,\W_G,\,\wh{\W}_G$ we also define equivariant 
cohomology algebras $\wh{H}_G(B),\,\H_G(B),\,\wh{\H}_G(B)$, 
which are modules over the rings $\ca{E}'(\g^*)^G,\,U(\g)^G,\,
\ca{E}'(G)^G$, respectively. 
 
The quantization map $\ca{Q}$ induces a vector space isomorphism 
$\ca{Q}:\,H_G(B)\to \ca{H}_G(B)$ and a linear map
$\ca{Q}:\,\wh{H}_G(B)\to \wh{\H}_G(B)$. Our second main result is 
that both of these maps are ring homomorphisms. 
For $B=\R$ this is the Duflo theorem for compact $G$.

In the case $B=\Om(M)$ where $M$ is a compact oriented $G$-manifold, there 
is a push-forward map $\int:\widehat{\H}_G(\Om(M))\to \ca{E}'(G)^G$.
In a sequel to this paper,  
we will explain how the localization formula 
in equivariant cohomology (see \cite{be:ze}, \cite{at:mom}) 
carries over to give a formula for the Fourier coefficients of 
this map.

In another sequel \cite{al:du} we use 
the equivariant cohomology groups $\widehat{\ca{H}}_G(\Om(M))$ 
to construct Liouville forms for Lie group valued 
moment maps, with applications to moduli spaces of 
flat connections on surfaces.
\vskip.2in

\noindent{\bf Acknowledgments.}  We would like to thank N. Berline,
R. Bott, P. Etingof, V. Guillemin, B. Kostant, A. Recknagel,
S. Sternberg, M. Vergne and C. Woodward for many stimulating discussions. 
We are particularly indebted to Michele Vergne for the important  
suggestion of using the duality between the Weil algebras 
$\wh{W}_G,\wh{\W}_G$ and spaces of 
differential forms $\Om(\g),\Om(G)$.
In particular, this idea improved the result about the 
ring structures in Section \ref{Kernels}.

\section{Review of the Weil algebra}
In this section we review H. Cartan's algebraic version of 
equivariant de Rham theory. Following more recent references
(e.g. \cite{ma:th,ka:br,gu:su})  we will phrase his theory 
in ``super''-terminology. Throughout this paper we work over
the field $\R$ of real numbers, and ``algebra'' will always mean  
algebra with a unit element.  
\subsection{Super-notation}
A super-vector space is a vector space $B$ together with a linear map
(parity operator) 
$\eps\in\End(B)$ satisfying $\eps^2=\on{Id}_B$. The $+1$ eigenspace
of $\eps$ is denoted $B_{even}$ and the $-1$ eigenspace $B_{odd}$. 
A homomorphism of super-vector spaces is a linear map preserving 
$\Z_2$-gradings.
Ungraded vector spaces $B$ will be viewed as ``super'' by putting 
$B=B_{even}$. The tensor product of two super-vector spaces $B_1,B_2$ is 
a super-vector space. If $B_j$ are super-algebras this tensor product 
denotes the super-($\Z_2$-graded) tensor product. Similarly 
commutators $[a,b]=\ad(a)b$ in a super-algebra
always mean super-commutators, derivations always mean
super-derivations.  The space $\Der(B)$ of derivations of a
super-algebra $B$ is a super-Lie algebra, and there is a super-Lie
homomorphism $\ad:\,B\to \Der(B)$.  An {\em odd} derivation $\d\in\Der(B)$ is
called a differential if $\d\circ \d=0$.

For any element $b\in B$ we denote by $b^L$ the operator of 
left multiplication by $b$ and by $b^R$ the operator by graded
right multiplication. Thus $\ad(b)=b^L-b^R$. For commutative 
super-algebras we drop the superscripts.

By a {\em graded}  super-vector space 
we mean a $\Z$-graded vector space $B^\star=\bigoplus_k B^k$, viewed 
as a super-vector space by putting $\eps=(-1)^k$ on $B^k$. A 
{\em filtered} super-vector space is defined to be 
a super-vector space $B=B^{(\star)}$, 
together with an increasing $\Z$-filtration 
$$\ldots\subseteq B^{(k)}\subseteq B^{(k+1)}\subseteq\ldots$$ 
such that 
$\eps (B^{(k)})\subseteq B^{(k)}$ for all $k$ and the 
associated graded space $\on{Gr}^*(B)$, with $\Z_2$-grading 
induced from $\eps$, is a graded super-vector space.

%

\subsection{G-differential algebras}
Let $G$ be a compact, connected Lie group with Lie algebra $\g$. 
Choose a basis $e_a$ of $\g$, and let
$e^a\in\g^*$ be the dual basis.  Let $f^c_{ab}$ be the structure
constants defined by
$$ [e_a,e_b]=f^c_{ab} e_c. $$
Let $\hat{\g}^\star$ be the graded super-Lie algebra defined as follows. 
As vector spaces 
$\hat{\g}^0\cong \hat{\g}^{-1}\cong \g$, while $\hat{\g}^1\cong\R$
and $\hat{\g}^k=\{0\}$ for $k\not=-1,0,1$. 
Letting $L_a,\iota_a$ denote the 
basis elements in $\hat{\g}^0,\,\hat{\g}^{-1}$ corresponding to $e_a$ 
and $\d$ the generator of $\hat{\g}^1$ the bracket relations are defined 
as
\beq
[L_a,\iota_b]&=&f^c_{ab}\iota_c,\\
{[ L_a, L_b ]} &= &f^c_{ab}L_c,\\
{[\iota_a,\d]}&=& L_a
\eeq
For any $G$-manifold $M$ there is a natural representation of $\hat{\g}$
on the commutative super-algebra of differential forms 
$B=\Om^\star(B)$, where $\iota_a$ are 
interpreted as contractions with generating vector fields
$$(e_a)_M:=\f{\p}{\p t}|_{t=0}\exp(-t e_a)^*$$
on $M$, $L_a$ as 
Lie derivatives, and $\d$ as the exterior differential. 
More generally, one defines: 
\begin{definition}
A {\em $G$-differential space}
is a super-space $B$, together 
with a super Lie algebra homomorphism $\rho:\,\hat{\g}\to \on{End}(B)$. 
The horizontal subspace 
$B_{hor}$ is the space fixed by $\hat{\g}^{-1}$, 
the invariant subspace $B^G$ is the 
space fixed by $\hat{\g}^0$, and the space $B_{basic}$ of basic 
elements is their intersection. 
If $B$ is graded/filtered and the map $\rho$ preserves 
degrees/filtration degrees we 
call $B$ a {\em graded/filtered $G$-differential space}. 
A $G$-differential algebra 
is a super-algebra $B$, together with a structure of a 
$G$-differential space such that $\rho$ takes values in 
derivations of $B$. 
\end{definition}

\begin{definition}
A homomorphism between (graded/filtered)  
$G$-differential 
spa\-ces/\-al\-ge\-bras $(B_1,\rho_1)$ and 
$(B_2,\rho_2)$ is a homomorphism of (graded/filtered)  
super-spaces/al\-ge\-bras
$\phi:\,B_1\to B_2$ with $A\circ \rho_1=\rho_2$.  
\end{definition}

The idea of a $G$-differential algebra was introduced by H. Cartan in \cite{ca:no}, 
Section 6 and  appears in the literature under various names. We 
follow the terminology from \cite{ku:eq}, although in this reference 
it is also required that $B$ carries a (Frechet) topology and 
the representation of $\hat{\g}^0=\g$ exponentiates to a differentiable 
$G$-action on $B$.
Obviously, the associated graded algebra of any filtered $G$-differential algebra
is a graded $G$-differential algebra. 

The basic examples for $G$-differential algebras are as follows. We will 
always write $\iota_a,L_a,\d$ in place of $\rho(\iota_a),\rho(L_a)$ and 
$\rho(\d)$.
\begin{example}[Trivial $G$-differential algebra]
The trivial $G$-differential algebra is $\R$ with the trivial representation 
$\rho=0$. If $(B,\rho)$ is any (possibly graded or filtered) 
$G$-differential algebra, the injection 
$$\R\to B,\ \ t\mapsto t\,I, $$ 
where $I\in B$ is the unit element,  
is a homomorphism of (graded/filtered) $G$-differential algebras. This follows 
since every derivation annihilates the unit element. 
\end{example}

\begin{example}[Exterior algebra]\label{Exterior}
Take $B=\wedge\g^*$, equipped with the coadjoint $G$-action.  Let
$y^a\in \wedge^1\g^*$ be the basis elements corresponding to $e^a$,
let $\iota_a=\iota(e_a)$ be the contraction,
$L_a$ the Lie derivative $L_a=-f_{ab}^c y^b \iota_c$, and let 
$\d$ be given by Koszul's formula
\begin{equation}\label{KoszulFormula}
\d=-\hh f^a_{bc} y^b y^c \iota_a=\hh y^a L_a. 
\end{equation}
Then $\d\circ\d=0$ and $[\iota_a,\d]=L_a$ so that 
$\wedge\g^*$ is a $G$-differential algebra. Since $G$ is compact,  
every cocycle has a representative 
in the invariant part $(\wedge^\star\g^*)^G$. Since 
$\d$ vanishes on this subspace it follows that 
$H^\star(\wedge\g^*,\d)\cong (\wedge^\star\g^*)^G$ . 
\end{example}

\begin{example}[Weil algebra]
Let $S\g^*$ be the symmetric algebra, equipped with the 
$G$-action induced by the coadjoint action of $G$ on $\g^*$. 
If we denote by $v^a\in S^1\g^*$ the generators corresponding to 
the basis $e^a\in\g^*$, the Lie derivative is given by the formula 
$$ L_a  v^c = - f_{ab}^c v^b.$$
The Weil algebra is the graded commutative super-algebra
$W_G=\oplus_{l=0}^\infty W^l_G$, where
$$ W^l_G=\bigoplus_{j+2k=l} S^k\g^* \otimes    \wedge^j\g^*.$$
We will write $v^a$ in place of $v^a\otimes 1 \in W^2_G$ and $y^a$ in place 
of $1\otimes y^a\in W^1_G$. 
Let $W^\star_G$ be equipped with the diagonal $G$-action and corresponding 
Lie derivative
$L_a=1\otimes L_a+L_a\otimes 1 \in\on{Der}^0(W_G)$. 
Let 
$\iota_a=1\otimes \iota_a\in\on{Der}^{-1}(W_G)$, and let 
the differential on $W_G$ be given by the formula
\begin{equation}\label{formula} 
\d =y^a (L_a\otimes 1) + (v^a -\hh f^a_{bc} y^b y^c) \iota_a.
\end{equation} 
On generators, 
\beq
\d y^a&=&v^a-\hh f^a_{jk}  y^j y^k,\\  
\d v^a&=&-f^a_{jk} y^j v^k.
\eeq
It is easily verified that 
$[\iota_a,\d]=L_a$ and $\d\circ\d =0$ 
so that $W^\star_G$ has the structure of a $G$-differential algebra.
Cartan shows in \cite{ca:la} that $(W_G,\d)$ is acyclic, 
that is, $H^k(W_G)=0$ for $k\not=0$, $H^0(W_G)=\R$.
The horizontal algebra is $(W_G)_{hor}\cong S\g^*$ and the 
basic subalgebra is the algebra of invariant polynomials
$(S\g^*)^G$. 
\end{example}
\subsection{Algebraic connections}
A $G$-differential algebra $B$ is called locally free 
if there exists an element 
$\theta=\theta^a e_a\in (B^{odd}\otimes\g)^G$, called (algebraic) 
{\em connection form} with 
$$\iota_a\theta=e_a.$$
If $B$ is graded/filtered we require in addition that $\theta\in 
(B^1\otimes\g)^G$ resp. $\theta\in (B^{(1)}\otimes \g)^G$.   
The invariance condition means that 
$ L_a \theta^c =-f_{ab}^c \theta^b.$
The $G$-differential algebra $B=\Om^\star(M)$ is locally free if 
and only if the
$G$-action on $M$ is locally free. In this case $\theta^a$ 
are connection forms in the usual sense. Note that 
the basic subcomplex $B_{basic}$ is naturally isomorphic to the 
differential forms on the quotient space:
$$ B^\star_{basic}= \Om^\star(M/G).$$
The exterior algebra $B=\wedge\g^*$ is locally free, with connection forms
$\theta^a=y^a$. The Weil algebra $W_G$ is locally free, with connection  forms
$\theta^a=y^a$. Together with acyclicity this 
motivates the interpretation of $W_G$ as the algebraic model for 
the de Rham complex of the classifying bundle $EG$.

\subsection{Equivariant cohomology}
Given any two $G$-differential algebras $B_1$ and $B_2$, the  
tensor product $B_1\otimes B_2$ is naturally a $G$-differential algebra.
If $B_1$ is locally free then the tensor product $B_1\otimes B_2$ 
is locally free. In particular, $W_G\otimes B$ is locally free 
for any $G$-differential algebra $B$. This motivates the following definition. 
\begin{definition}
Let $B$ be any $G$-differential algebra. The {\em equivariant cohomology} 
of $B$ is super-algebra defined as the 
cohomology of the basic subcomplex,
$$ H_G(B):=H((W_G\otimes B)_{basic}).$$ 
\end{definition}
If $B$ is a graded/filtered super-algebra then $H_G(B)$ inherits 
a grading/filtration.

In the special case $B^\star=\Om^\star(M)$, and if $G$ is compact, 
it can be shown that $H^\star_G(\Om^\star(M))$ equals the 
topological equivariant cohomology $H^\star_G(M):=H^*(EG\times_G M)$. 
For non-compact $G$ this statement is false in general. 
\begin{remarks}
\begin{enumerate}
\item
If $B$ is locally free, $H_G(B)=H(B_{basic})$.
\item
The equivariant cohomology of the trivial $G$-differential algebra is
$H_G(\R)=(S\g^*)^G$.  
\end{enumerate}
\end{remarks}

Any homomorphism of $G$-differential spaces/algebras $\phi:\,B_1\to B_2$ induces a 
super-space/algebra homomorphism $H_G(B_1)\to H_G(B_2)$. It follows that 
for any $G$-differential algebra $B$, 
the natural homomorphism $\R\to B$ induces an algebra homomorphism 
$H_G(\R)\to H_G(B)$ 
making $H_G(B)$ into a module over the ring of invariant polynomials
$H_G(\R)=(S\g^*)^G$. If $M$ is an 
oriented compact manifold, the integration map $\int_M:\Om^\star(M)\to \R$ 
is a chain map, and therefore defines an equivariant linear map
$\int_M:\,H_G(M)\to S(\g^*)^G$.

\begin{definition} \label{ChainHomotopic}
Two homomorphisms $\phi_1,\phi_2:\,B_1\to B_2$  
between $G$-differential spaces 
are called $G$-chain homotopic if there an odd linear map
$h:\,B_1\to B_2$ which commutes with contractions 
$\iota_a$ and Lie derivatives $L_a$ and which satisfies 
$[\d,h]=\phi_1-\phi_2$. 
\end{definition}
If $\phi_1,\phi_2:\,B_1\to B_2$ are chain homotopic homomorphisms 
of $G$-differential spaces, then the induced maps in cohomology 
coincide.

\section{The non-commutative Weil algebra}
All of the examples of $G$-differential algebras discussed in the previous section 
are super-commutative. In this section we will consider two non-commutative 
examples: The Clifford algebra $\Cl(\g)$ and its tensor product 
with the universal enveloping algebra, $U(\g)\otimes \Cl(\g)$. 
\subsection{The Clifford algebra}
Suppose that $\g$ comes equipped with an invariant inner 
product, used to identify $\g\cong\g^*$. 
Let $\Cl(\g)$ be the Clifford algebra of $\g$, i.e. the quotient of
the tensor algebra by the ideal generated by all
$(\mu,\mu)-2\mu\otimes\mu$ with $\mu\in\g$. It inherits from the tensor
algebra a natural $\Z_2$-grading and filtration, 
$$\R=\Cl^{(0)}(\g)\subset \Cl^{(1)}(\g)\subset\ldots.$$
The filtration and grading are 
compatible, so that $\Cl(\g)$ is a filtered super-algebra.
There is a canonical isomorphism 
$$\Gr^\star(\Cl(\g))\cong \wedge^\star\g.$$  
Let us choose 
the basis $e_a$ of $\g$ to be orthonormal so that $e_a=e^a$. Using the 
metric to pull down indices we write $f_{abc}=f^c_{ab}$. 
If we denote the generators of $\Cl(\g)$ corresponding to 
$e_a$ by $x_a$ the defining relations are 
$ [x_a,x_b]=\delta_{ab}.$  

The {\em symbol map} $\sig:\,\Cl(\g)\to \wedge\g$ is the vector 
space isomorphism given on generators by 
$$ \sig(x_{j_1}\ldots x_{j_k})=y_{j_1}\wedge\ldots \wedge y_{j_k}$$
for ${j_1}<\ldots<{j_k}$. It does not 
depend on the choice of basis. 

Henceforth we will sometimes drop $\sig$ from the notation and 
just think of $\Cl(\g)$ as $\wedge\g$ with a 
different product structure 
$\odot$, defined by 
$$ y\odot y'= \sig(\sig^{-1}(y)\sig^{-1}(y')).$$
The operators of Clifford multiplication by $x_a$ from the left 
or right become
\begin{equation}\label{LeftRight}
x_a^L=y_a+\hh \iota_a,\ \ x_a^R=y_a-\hh \iota_a.
\end{equation}
The relation between Clifford 
multiplication and exterior multiplication is described 
in \cite{ko:cl}, Theorem 16: 
Let $\iota_a^1$ and $\iota_a^2$ be the contraction operators for 
the first resp. second factor in $\wedge\g\otimes\wedge\g$. 
\begin{lemma}
\label{KostantLemma}
The algebra structure on $\wedge\g$ induced by the symbol map 
is the map $\odot:\wedge\g\otimes\wedge\g\to \wedge\g$ given by 
composition of the operator $\exp(-\hh \iota_a^1\iota_a^2)$ 
on $\wedge\g\otimes\wedge\g$ with exterior multiplication on 
$\wedge\g$. 
\end{lemma}
\begin{proof}
Since the product $\odot$ is associative, it suffices to check
$\sig(x\,x')=\sig(x)\odot\sig(x')$ for the case $x=x_j$ for some $j$, 
and $x'=x_{j_1}\ldots x_{j_k}$ with ${j_1}<\ldots<{j_k}$. In both 
of the sub-cases
$j\in\{j_1,\ldots ,j_k\}$ and $j\not\in\{j_1,\ldots ,j_k\}$ this 
is easily verified.
\end{proof}
Let us now consider the $G$-action on $\Cl(\g)$
induced by the adjoint action on 
$\g$.
The corresponding Lie derivative is a commutator $L_a=\ad(g_a)$ 
with 
\begin{equation}
\label{Ga}
g_a=-\hh f_{ars} x_r x_s\in \on{Cl}^{(2)}(\g).
\end{equation}
Let 
$$ \gamma=\f{1}{3}x_a g_a 
= -\f{1}{6}f_{abc}x_a x_b x_c \in \Cl^{(3)}(\g)^G.$$
\begin{proposition}
The Clifford algebra $\Cl(\g)$ with derivations 
\beq 
\iota_a&=&\ad(x_a),\\
L_a&=&\ad(g_a),\\ 
\d&=& \ad(\gamma)
\eeq
is a filtered $G$-differential algebra, having $\wedge \g$ 
as its associated graded $G$-differential algebra. The cohomology is 
trivial in all filtration degrees (except if $\g$ is abelian, in which case 
$\d=0$). 
\end{proposition}

\begin{proof}
The required commutation relations for the operators 
$L_a,\iota_a$ and $\d$ follow from those for the elements 
$g_a,x_a$ and $\gamma$. For example $[g_a,\gamma]=L_a\gamma=0$
shows $[L_a,\d]=0$, and $[x_a,\gamma]=\iota_a\gamma=g_a$ 
shows $[\iota_a,\d]=0$. The only non-trivial commutation relation 
to check is that $[\d,\d]=0$, or equivalently that
$[\gamma,\gamma]=2\gamma^2$ is in the center of $\Cl(\g)$. 
In fact $\gamma^2$ is a scalar (cf. Kostant \cite{ko:cl}): 
\begin{equation}\label{Square} 
\gamma^2=-\f{1}{48}\,f_{abc}f_{abc}.
\end{equation}
To see this we compute the symbol $\sig(\gamma^2)
=\sig(\gamma)\odot\sig(\gamma)$
using Lemma \ref{KostantLemma}. Here
$\sig(\gamma)=-\f{1}{6}f_{abc}y_a y_b y_c$. 
The terms $\sig(\gamma)$ and $\iota_a\iota_b\sig(\gamma)$ square to
zero since they have odd degree, and the term
$\iota_a\sig(\gamma)=\sig(g_a)$ squares to zero by the Jacobi
identity. Since
$$ \f{1}{3!}(-\hh \iota_a^1\iota_a^2)^3 \sig(\gamma)\otimes\sig(\gamma) 
=-\f{1}{48} (\iota_a^1\iota_a^2)^3 \sig(\gamma)\otimes\sig(\gamma) 
=-\f{1}{48} f_{abc}f_{abc}
$$
Equation \eqref{Square} follows. 
This shows that the operators $\d,L_a$ and $\iota_a$ define a 
representation of $\hat{\g}$ on $\Cl(\g)$.
To verify that these operators
have the required filtration degrees we check on generators:
$$ L_a x_b = f_{abc} x_c, \ \ \iota_a x_b = \delta_{ab}, \ \
     \d x_a=-\hh f_{ars} x_r x_s .
$$

These equations show also that the associated graded $G$-differential 
algebra
on $\wedge^\star\g=\on{Gr}^\star(\Cl(\g))$ is the standard one. 
Finally, to see that the cohomology of $(\Cl(\g),\d)$ is trivial 
(if $\g$ is non-abelian) we note that 
$$[\d,\gamma]=[\gamma,\gamma]=2\gamma^2=\f{-1}{24} f_{abc}f_{abc},$$
so that 
$$ H:=\f{-24}{f_{abc}f_{abc}}\  \gamma $$
is a homotopy operator. 
\end{proof}

Below we will need the following description of the differential 
$\d=\ad(\gamma)$ in terms of the identification $\sig:\,\Cl(\g)\cong
\wedge\g$. It shows that $\ad(\gamma)$ 
is the Koszul differential plus an extra cubic term:
\begin{proposition}\label{ClifDConj}
Under the identification $\Cl(\g)\cong\wedge\g$ by the symbol map, 
$$ \ad(\gamma) =   
-\hh f_{abc}
y_b y_c \iota_a - \f{1}{24} f_{abc} \iota_a\iota_b\iota_c.
$$
\end{proposition}
\begin{proof}
Using \eqref{LeftRight} 
we compute 
\beq 
\gamma^L&=&-\f{1}{6}f_{abc}\big((y_a+\hh \iota_a)
(y_b+\hh\iota_b)(y_c+\hh\iota_c)\big)\\
\gamma^R&=&-\f{1}{6}f_{abc}\big((y_a-\hh \iota_a)
(y_b-\hh \iota_b)(y_c-\hh\iota_c)\big).
\eeq
Taking the difference $\ad(\gamma)=\gamma^L-\gamma^R$, 
the terms cubic and linear in $y_a$'s cancel, and the remaining 
terms yield the Proposition.  
\end{proof}

\subsection{The non-commutative Weil algebra}
Let $U(\g)$ be the universal enveloping algebra, with its natural 
filtration $\R=U^{(0)}(\g)\subset U^{(1)}(\g)\subset\ldots$. 
The associated graded 
algebra is  $\Gr^\star U(\g)=S^\star\g$. We denote the generators of 
$U(\g)$ corresponding to $e_a\in\g$ by 
$u_a$. The extension to $U(\g)$ of the adjoint action 
has Lie-derivative $ L_a=\ad(u_a)$.
Define a filtered super-algebra $\W_G^{(\star)}$, where 
$$ \W^{(l)}_G=\bigoplus_{j+2k=l} U^{(k)}(\g)\otimes \Cl^{(j)}(\g),$$
equipped with the diagonal $G$-action.
Then $\Gr^\star(\W_G)=W_G^\star$ is just the usual Weil algebra, using 
the identification $\g^*\cong\g$. We will show that $\W_G$ has the 
structure of a filtered $G$-differential algebra.

The generators $u_a=u_a\otimes 1$ and $x_a=1\otimes x_a$ have 
degree $2$ and $1$, respectively.
The Lie derivatives are $ L_a=\ad(u_a+g_a)$, and the
derivations $\iota_a$ on $\Cl(\g)$ extend to derivations $\iota_a$ 
on $\W_G$. The basic subspace of $\W_G$ for the derivations $\iota_a$ and 
the $G$-action is $(\W_G)_{basic}=U(\g)^G\otimes 1$. In particular, 
$(\W_G)_{basic}$ is a subspace of the center $Z(\W_G)$.
Introduce an element $\D\in (\W_G^{(3)})^G$ by 
$$ \D:= x_a u_a+\gamma.$$ 
\begin{proposition}\label{SquareD}
The square $\D^2$  is  given by 
$$\D^2=\hh u_a u_a+\gamma^2
=\hh u_a u_a-\f{1}{48} f_{abc}f_{abc}.$$
\end{proposition}
\begin{proof}
We calculate:
\beq
\D^2=\hh [\D,\D]&=&
\hh [u_a x_a,  u_b x_b] + \gamma^2 + [u_a x_a ,\gamma]\\
&=& \hh u_a u_a + \hh [u_a,u_b]  x_a x_b+\gamma^2 
+u_a [x_a,\gamma]\\
&=&\hh u_a u_a + \hh f_{abc} u_c x_a x_b + \gamma^2
+ u_a g_a \\
&=& \hh u_a u_a + \gamma^2.
\eeq
\end{proof}

\begin{remark}
As suggested by the formula for $\D^2$, the element $\D$ may 
be viewed as a Dirac operator. Indeed let $\Cl(T^*G)$ be the 
Clifford algebra bundle for the left-invariant Riemannian 
metric on $G$ which coincides with the given inner
product on $\g=T_eG$. Identify $x_a$ with left-invariant 
sections of $\Cl(T^*G)$ and $e_a$ with left-invariant 
(co-)vector fields. Let $\nabla_a^L=\nabla_{e_a}^L$ denote covariant 
derivatives with respect to left-trivialization. The operators 
$\nabla_a^L$ satisfy the 
commutation relations $[\nabla_a^L,\nabla_b^L]=f_{abc}\nabla_c^L$
and $[\nabla_{e_a}^L,x_b]=0$. Hence setting $u_a=\nabla_a^L$ 
we identify $U(\g)\otimes \Cl(\g)$ with an algebra of 
left-invariant differential operators acting on the bundle
$\wedge T^*G$. Furthermore $\nabla_a=\nabla^L_a-\f{1}{6} f_{abc}x_bx_c$ 
becomes a left-invariant connection and $\D=x_a \nabla_a$ the 
associated Dirac operator. (It is not the de-Rham Dirac 
operator, however.) The above Dirac operator 
appears in the Physics literature, see
Fr\"ohlich, Grandjean and Recknagel \cite{fr:su} (Section 3).

A generalization of this Dirac operator, with fascinating 
applications to representation theory was considered by Kostant 
\cite{ko:cu}. 
\end{remark}
Since $\D^2\in (\W_G)_{basic}$ is a central element, 
$\d:=\ad(\D)$ is a differential. On generators, 
\beq
\d x_a&=& u_a - \hh f_{ajk} x_j x_k,\\
\d u_a&=& -f_{ajk} x_j u_k,
\eeq
which shows that $\d$ has filtration degree 1.

\begin{theorem}
The non-commutative Weil algebra  $\W_G$, equipped with derivations 
\beq
\iota_a&=&\ad(x_a),\\
L_a&=&\ad(u_a+g_a),\\
\d&=&\ad(\D)
\eeq
is a filtered $G$-differential algebra, with associated graded $G$-differential 
algebra equal to the standard Weil algebra $W_G$ (identifying $\g^*\cong\g$). 
\end{theorem}

\begin{proof}
We have already shown that $[\d,\d]=0$. 
Since $\D$ is $G$-invariant, $[g_a,\D]=L_a\D=0$ so that 
$[L_a,\d]=0$. Furthermore 
$$ [x_a,\D]= u_a-\hh f_{ars} x_r x_s =u_a+g_a$$
shows $[\iota_a,\d]=L_a$. All other commutators are obvious.
The formulas for $\d,\iota_a,L_a$ on generators show that 
the associated graded $G$-differential algebra is just the standard one on 
$W_G$.
\end{proof}  
To compare the differential on $\W_G$ with the differential on the 
commutative Weil algebra $W_G$, use the symbol map 
to identify $\W_G=U(\g)\otimes\Cl(\g)\cong U(\g)\otimes \wedge\g$. 
\begin{proposition}\label{ConjSigma}
Under the identification $\sig:\,
\W_G\cong U(\g)\otimes \wedge\g$, 
the non-commutative Weil differential
is given by the formula 
$$
\d^\W = y_a (L_a\otimes 1)+  
\Big(\f{u_a^L+u_a^R}{2} -\hh f_{abc}
y_b y_c\Big) \iota_a - \f{1}{24} f_{abc} \iota_a\iota_b\iota_c.
$$
\end{proposition}
\begin{proof}
Using \eqref{LeftRight} we compute: 
\beq
\ad(x_a u_a)
&=&(y_a+\hh \iota_a)u_a^L-(y_a-\hh \iota_a)u_a^R\\
&=&(u_a^L-u_a^R)y^a+\hh(u_a^L+u_a^R)\iota_a.
\eeq
Together with Proposition \ref{ClifDConj} this proves
Proposition \ref{ConjSigma}.
\end{proof}

\begin{definition}
If $B$ is any $G$-differential algebra
define the cohomology $\H_G(B)$ 
as the super-algebra
$$ \H_G(B):=H((\W_G\otimes B)_{basic},\d).$$
If $B$ carries a filtration $B=B^{(\star)}$ then 
$\H_G(B)$ inherits a filtration $\H^{(\star)}(B)$. 
\end{definition}
Any homomorphism $B_1 \to B_2$ of (filtered) $G$-differential algebras
induces an homomorphism of (filtered) super-algebras
$\H_G(B_1)\to \H_G(B_2)$. In particular, for any 
$G$-differential algebra $B$ the natural 
embedding $\R\to B^{(0)},\ t\mapsto tI$ 
makes $\H_G(B)$ into an algebra over the 
ring of Casimir operators $\H_G(\R)=U(\g)^G$. If $M$ is an 
oriented compact manifold, the integration map $\int_M:\Om(M)\to \R$ 
induces a linear map in cohomology, 
$\int_M:\,\H_G(\Om (M))\to U(\g)^G$.

\section{The Cartan model}\label{CartanModel}
In \cite{ca:la}, Section 6 Cartan shows that
the equivariant cohomology of a $G$-differential algebra can be computed in an 
equivalent, simpler ``model'' $C_G(B):=(S\g^*\otimes B)^G$ 
with differential $\d_G=1\otimes\d-v^a\otimes \iota_a$. 
The goal of this section is 
to develop a similar Cartan model for $\ca{H}_G(B)$. 
We first review the standard Cartan model, following the exposition in 
Guillemin-Sternberg \cite{gu:su}. 
\subsection{Cartan model for $H^\star_G(B)$.}
A simple and transparent way of obtaining the Cartan model from  
the Weil model was found by Kalkman \cite{ka:br}, following earlier work 
by Mathai-Quillen \cite{ma:th}. Let the Kalkman operator $\phi$ be 
defined by
$$\phi=\exp(y^a\otimes\iota_a):\, W_G\otimes B \to  W_G\otimes B .$$
One checks that $y^a\otimes\iota_a\in\on{Der}(W_G\otimes B)$ 
is an even derivation so that $\phi$ defines an algebra isomorphism
(preserving degrees/filtration degrees if $B$ is graded/filtered).
The key property of the Kalkman map $\phi$ is as follows.
\begin{proposition}[\cite{ka:br}, Eq. (1.23)] 
\label{MathaiQuillen}
The conjugate of the contraction operator $\iota_a=
\iota_a\otimes 1+1\otimes\iota_a$ under the map $\phi$ is
$$ \Ad_{\phi}\iota_a=\iota_a\otimes 1 .$$
\end{proposition}
\begin{proof}
This follows by writing $\Ad_\phi=\sum_j \f{1}{j!} \ad^j(y^r\otimes\iota_r)$ 
and 
\beq 
\ad(y^r\otimes\iota_r)(\iota_a\otimes 1+1\otimes \iota_a)&=&
-1\otimes\iota_a,\\
\f{1}{2!}
\ad(y^r\otimes\iota_r)(\iota_a\otimes 1+1\otimes \iota_a)&=&
0.
\eeq
\end{proof}
Since the kernel of the operators $\iota_a\otimes 1 $ is just
$S\g^*\otimes B$, it follows that $\phi$ restricts to an algebra
isomorphism $(W_G \otimes B)_{hor}\to S\g^*\otimes B$, and by
equivariance it restricts further to an algebra isomorphism 
$(W_G\otimes B)_{basic}\to C_G(B)=(S\g^*\otimes
B)^G$.

The algebra isomorphism $(W_G\otimes B)_{basic}\to C_G(B)$ 
is due to Cartan \cite{ca:no} and 
Mathai-Quillen \cite{ma:th} while the extension to an algebra 
automorphism of $W_G\otimes \Om(M)$ is due to Kalkman 
\cite{ka:br}. It can also be described as follows. Let 
$$ P_{hor}=\prod_\beta \iota_\beta y^\beta:\,W_G\to 
S(\g)\hra W_G$$
be the horizontal projection. On 
the invariant subspace $(W_G\otimes B)^G$, the Kalkman map 
agrees with the operator $P_{hor}\otimes 1$ 
since
$$
\exp(y^a\otimes \iota_a)=
\prod_\beta(1+y^\beta \otimes \iota_\beta)=\prod_\beta(1-
y^\beta\iota_\beta \otimes 1)=\prod_\beta(\iota_\beta y^\beta
\otimes 1)=P_{hor}\otimes 1.
$$
(We use the convention that we sum over roman indices but not over 
greek indices.) 
The Cartan differential $\d_G$ is defined by the condition 
$$ \d_G\circ (P_{hor}\otimes 1)=
(P_{hor}\otimes 1) \circ (\d\otimes 1+1\otimes\d)$$
on $(W_G\otimes B)_{basic}$. 
Recall formula \eqref{formula} for 
$\d\otimes 1$. Application of $P_{hor}\otimes 1$ annihilates all terms
involving $y_a's$, and on the subspace 
$(W_G\otimes B)_{basic}$, the operator $v^a(\iota_a\otimes 1)$
may be replaced by $-v^a(1\otimes\iota_a)$ which then commutes 
with $(P_{hor}\otimes 1)$. Hence 
$$ (P_{hor}\otimes 1) \circ (\d\otimes 1+1\otimes\d)
=(1\otimes \d -v^a \otimes \iota_a)\circ (P_{hor}\otimes 1)
$$ 
on $(W_G\otimes B)_{basic}$. Therefore
$$ \d_G=1\otimes \d -v^a \otimes \iota_a:\,C_G(M)\to C_G(M).$$
Thus, the equivariant cohomology $H_G(B)$ can be computed as the 
cohomology of the complex $(C_G(B),\d_G)$: 
$$H_G(B)=H(C_G(B),\d_G).$$

\subsection{Cartan model for $\ca{H}_G(B)$}
We now describe the Cartan model for the equivariant cohomology 
$\ca{H}_G(B)$. Let $B$ be any $G$-differential algebra. The non-commutative 
analogue of the Kalkman isomorphism is the $G$-equivariant operator
$$\phi:=\exp(x_a\otimes \iota_a):\,\W_G\otimes B\to \W_G\otimes B.$$
Note that if $B$ is filtered, the map $\phi$ preserves filtration 
degrees. %
By a calculation parallel to that for the usual Weil algebra, 
$\phi$ coincides  on $(\W_G\otimes B)_{hor}$ with the operator
$P_{hor}\otimes 1$ where $P_{hor}:\,\W_G\to U(\g)$ is the 
horizontal projection
$$P_{hor}=\prod_\beta\,\iota_\beta\, x_\beta^L
=\prod_\beta\,\iota_\beta\, y_\beta
.$$ 
Likewise, the proof of 
Proposition \ref{MathaiQuillen} carries over to show that
$$\Ad_\phi \iota_a=\iota_a\otimes 1.$$
Thus $\phi$ defines a vector space isomorphism 
$(\W_G\otimes B)_{hor}\to U(\g)\otimes B$, and by equivariance
$$(\W_G\otimes B)_{basic}\cong (U(\g)\otimes B)^G.$$ 
We call $\ca{C}_G(B):=(U(\g)\otimes B)^G$ the non-commutative
Cartan model. The new differential $\d_G$ on $\ca{C}_G(B)$ is once 
again be obtained from the condition 
$$ \d_G\circ (P_{hor}\otimes 1)=
(P_{hor}\otimes 1) \circ (\d\otimes 1+1\otimes\d)$$
on $(\W_G\otimes B)_{basic}$. One can read off an expression 
for $\d_G$ from Proposition \ref{ConjSigma}: Application of
$(P_{hor}\otimes 1)$ annihilates all terms involving $y_a$'s, 
and on the subspace $(\W_G\otimes B)_{basic}$ we may replace
$\iota_a\otimes 1$ by $-1\otimes\iota_a$ which then commutes 
with $(P_{hor}\otimes 1)$. 

Hence we have shown:
\begin{proposition}
The differential $\d_G$ on $\ca{C}_G(B)=(U(\g)\otimes B)^G$ 
is given by 
$$ \d_G=1\otimes \d -\hh (u_a^L+u_a^R)\otimes
\iota_a+\f{1}{24} f_{abc} (1\otimes \iota_a \iota_b\iota_c).
$$
In particular, $\d_G\circ \d_G=0$ on $\ca{C}_G(B)$.
\end{proposition}
It follows that the equivariant cohomology $\ca{H}_G(B)$ 
of any $G$-differential algebras is equal to the cohomology of the 
complex $\ca{C}_G(B)$ with differential $\d_G$. 

\subsection{Ring structure of the non-commutative Cartan model}
In contrast to the commutative Weil model the non-commutative
Kalkman map is not an algebra isomorphism (since 
$x_a\otimes\iota_a $ is not a derivation). In this section we 
compute the new
algebra  structure on the Cartan model $\ca{C}_G(B)$ induced by the 
Kalkman map.

Suppose $B$ is a $G\times G$-differential algebra. Passing to 
the diagonal action it becomes a $G$-differential algebra. 
The ring structure for $W_G$ defines a natural map 
$(W_{G\times G}\otimes B)_{basic}\to (W_G\otimes B)_{basic}$, 
where on the left hand side mean basic elements for the $G\times G$-action 
and on the right hand side for the $G$-action. Correspondingly, 
we obtain a chain map between Cartan models, 
$$ \gamma:\,\ca{C}_{G\times G}(B)\to \ca{C}_G(B).$$

\begin{proposition}\label{Ring}
For any $G\times G$-differential algebra $B$, the 
map between Cartan models induced by the 
ring structure of $\W_G$ reads
$$  (\on{Mult}_{U(\g)}\otimes 1)\circ 
\exp(-\hh (1\otimes \iota_a^1\iota_a^2)):\,\ca{C}_{G\times G}(B)\to 
\ca{C}_G(B)$$
where $(\on{Mult}_{U(\g)}$ is the multiplication map for the universal 
enveloping algebra $U(\g)$.
\end{proposition}
\begin{proof}
Let $\Mult_{\W_G}=\Mult_{U(\g)}\otimes \Mult_{\Cl(\g)}$ be
the multiplication map for $\W_G$. 
The map $\gamma$ is defined by the condition $\gamma\circ 
(P_{hor}\otimes 1)=(P_{hor}\otimes 1)\circ (\Mult_{\W_G}\otimes 1)$ 
on basic elements $(\W_{G\times G}\otimes B)_{basic}$. 
According to Lemma \ref{KostantLemma}, 
under the identification $\sig:\,\Cl(\g)\cong \wedge\g$, 
Clifford multiplication is given by composition of 
the operator $\exp(-\hh \iota_a^1\iota_a^2)$ followed by 
wedge product. On basic elements, we can replace $\iota_a^j\otimes 1$
by $-(1\otimes \iota_a^j)$, and therefore 
$\exp(-\hh (\iota_a^1\iota_a^2\otimes 1))$ with 
$\exp(-\hh (1\otimes \iota_a^1\iota_a^2))$, which then commutes 
with $P_{hor}\otimes 1$.  
\end{proof}

For any $G$-differential algebra $B$ define an associative ring 
structure on $\ca{C}_G(B)$ as follows. 
Let $\iota_a^1,\iota_a^2$ be the contraction 
operators on $B\otimes B$ with respect to the first and second 
$G$-factor. Define a new ring structure $\odot$ on 
$B$ as a composition of the operator $\exp(-\hh\iota_a^1\iota_a^2)$ 
on $B\otimes B$ with  multiplication map $B\otimes B\to B$. 
The ring structure extends to a ring structure $\odot$ on the Cartan model
$\ca{C}_G(B)$. It follows from Proposition \ref{Ring} (applied to 
$B\otimes B$) that the ring structure $\odot$ on $\ca{C}_G(B)$ 
is a chain map for $\d_G$, or equivalently that 
$\d_G$ is a derivation for $\odot$.

\section{Generalized coefficients}
We will often find it useful to embed commutative and
non-commutative Weil algebras into 
``Weil algebras  with generalized coefficients''. 
This point of view will make many of the subsequent 
constructions much more natural, and also it is  
dictated by the applications to moment map theory developed 
in \cite{al:du}. 

We point out that already for the 
commutative Weil model, we introduce 
generalized coefficients in a way different from Kumar-Vergne
\cite{ku:eq} since we are working in the Fourier transformed 
picture.

\subsection{The Weil algebra with generalized coefficients}
Let $\ca{E}'(\g^*)$ be the convolution algebra of compactly supported 
distributions on $\g^*$. The subalgebra of distributions with 
support at the origin is canonically isomorphic to 
$S\g^*$, by identifying the generators $v^a$ with
$$ v^a= \f{\d}{\d t}\Big|_{t=0} \delta_{t e^a}=-\f{\p}{\p \mu_a} \delta_0.$$ 
%
We view the coordinate functions $\mu_a$ as multiplication 
operators $\l\mu_a u,\phi\r=\l u,\mu_a\phi\r$. Then 
$ [\mu_a,v^b]=\delta_a^b$ and   
the Lie derivative $L_a$ for the conjugation action on $\ca{E}'(\g^*)$
can be written
\begin{equation}\label{LieDerivative}
L_a=f_{ab}^c v^b\mu_c.
\end{equation}

The {\em Weil algebra with generalized coefficients} is the 
$G$-differential algebra
$$ \widehat{W}_G=\ca{E}'(\g^*)\otimes \wedge\g^*,$$
where the differential is given by formula \eqref{formula}. 
The multiplication map on $\widehat{W}_G$ extends to a map 
$$ \on{Mult}_W:\widehat{W}_{G\times G}\to \widehat{W}_G$$ 
from the completion $\widehat{W}_{G\times G}$ of $\widehat{W}_G\otimes
\widehat{W}_G$. This map is the direct product of the push-forward under the 
addition map $\on{Add}_*:\,\ca{E}'(\g^*\times\g^*)\to
\ca{E}'(\g^*)$ on the distribution part
and the multiplication map on the exterior algebra part.

Given a $G$-differential algebra $B$ we define the cohomology 
with generalized coefficients as the super-algebra
$$\widehat{H}_G(B):=H((\widehat{W}_G\otimes B)_{basic}).$$
It is naturally a module over the ring $ \widehat{H}_G(\R)=
\ca{E}'(\g^*)^G$ of invariant compactly supported 
distributions. The embedding $S(\g^*) \hra \ca{E}'(\g^*)$ induces
a natural embedding $W_G \hra \widehat{W}_G$ and
for any $G$-differential algebra $B$
a map in cohomology $H_G(B) \rightarrow \widehat{H}_G(B)$.

The Weil differential simplifies if we conjugate it by the function
$\tau_0\in C^\infty(\g^*,\wedge\g^*)$ given as
\begin{equation}\label{DefTau0}
 \tau_0(\mu)=\exp(-\hh f^a_{bc}\,\mu_a y^b y^c) 
\end{equation}
where $\mu=\mu_a e^a$. 
\begin{lemma}\label{Conjugate}
The conjugate of the differential $\d$ and the contraction $\iota_a$ 
by the function $\tau_0$ are given by 
\beq \Ad(\tau_0^{-1})\,\d&=&v^a \iota_a\\
\Ad(\tau_0^{-1})\,\iota_a&=&\iota_a-f_{ab}^c\mu_c y^b 
\eeq
\end{lemma}

\begin{proof}
We compute $\Ad(\tau_0)(v^r \iota_r)=\sum_j \f{1}{j!}
\ad^j(-\hh f^a_{bc}\,\mu_a y^b y^c)(v^r \iota_r)
$: 
\beq 
\ad(-\hh f^a_{bc}\,\mu_a y^b y^c)(v^r \iota_r)&=&
-\hh f^a_{bc}\,y^b y^c \iota_a + f^a_{bc}\,\mu_a v^b y^c \\
&=& -\hh f^a_{bc}\,y^b y^c \iota_a   +(L_c\otimes 1) y^c ,\\
\f{1}{2!}\ad^2(-\hh f^a_{bc}\,\mu_a y^b y^c)(v^r \iota_r)&=&
\hh f^a_{bc} f^c_{rs}\,\mu_a y^b y^r y^s=0
\eeq
where the last term vanishes by the Jacobi identity for $\g$.
Similarly 
\beq \ad(\hh f^c_{bd} \mu_c y^b y^d)\iota_a&=&
f^c_{ba} \mu_c y^b\\
\f{1}{2!}\ad^2(\hh f^c_{bd} \mu_c y^b y^d)\iota_a&=&0.
\eeq  
\end{proof}

Lemma \ref{Conjugate} has the following interpretation. 
Consider the natural pairing of $\wh{W}_G$ with the space
$\Om(\g^*)=C^\infty(\g^*)\otimes\wedge\g$ of smooth differential 
forms on $\g^*$. All the operators on $\wh{W}_G$ that we are
interested in are dual to certain operators on $\Om^\star(\g^*)$. 
In particular, 
$$ 
v^a=\big(\f{\p}{\p \mu_a}\big)^*,\ \ 
y^a=-\iota\big(\f{\p}{\p\mu_a}\big)^*,\ \ 
\iota_a=-(\d\mu_a)^*,\ \ 
\mu_a=(\mu_a)^*.
$$
and by Lemma \ref{Conjugate}, 
\begin{eqnarray}
\Ad(\tau_0^{-1})\d&=&v^a\iota_a=-\d^*\label{A}\\
\Ad(\tau_0^{-1})\iota_a&=&\iota_a-f_{ab}^c\mu_cy^b=
-(\d\mu_a+\iota_a)^*.\label{B}
\end{eqnarray}
Note that $\Om^\star(\g)$ with derivations
$\ti{\iota}_a=\iota_a+\d\mu_a,L_a,\d$ is a $G$-differential space, so
that $\wh{W}_G$ is (in an obvious sense) the dual $G$-differential
space.  Letting $\Mult_W:\,\wh{W}_G\otimes\wh{W}_G\to \wh{W}_G$ denote
the multiplication in the Weil algebra, the composition
$\tau_0^{-1}\circ \Mult_W\circ (\tau_0\otimes\tau_0)$ is dual to the
co-product $\Om^*(\g^*)\to \Om(\g^*)\otimes\Om(\g^*)$ given by
pull-back under the addition map $\on{Add}:\,\g^*\times\g^*\to\g^*$.

As applications we have: 
\begin{proposition}
For all {\em fixed} $\mu$, the element 
$\tau_0(\mu) \delta_\mu\in \widehat{W}_G$ 
is closed.
\end{proposition}
\begin{proof}
We verify:
$$ \d \tau_0(\mu) \delta_\mu(\nu) =
\d \tau_0(\nu) \delta_\mu(\nu) =
\tau_0(\nu) (v^a\iota_a)\delta_\mu(\nu) =0.$$
\end{proof}

\begin{proposition}
The element
$$ \Lambda_0=\exp(-y_a\d\mu_a)\tau_0(\mu)\delta_\mu
\in \wh{W}_G\otimes\Om(\g^*)
$$ 
has the properties 
$$ \d \Lambda_0=0,\ \ \iota_a\Lambda_0=-\d\mu_a\Lambda_0.$$
\end{proposition}
\begin{proof}
Since 
$\exp(-y_a\d\mu_a)\delta_\mu\in\wh{W}_G\otimes\Om(\g^*)$ 
can be viewed as
the kernel of the identity map which obviously commutes
with contractions and with the differential,
this follows from Equations \eqref{A} and \eqref{B}.
\end{proof}

As another application of Lemma \ref{Conjugate}, let us show 
that $\widehat{W}_G$ is acyclic. 
The space of compactly supported distributions is a direct sum 
$$ \ca{E}'(\g^*)=\R\oplus \ca{E}'(\g^*)_+ $$
where $\R$ is embedded as multiples of $\delta_0$ and 
$\ca{E}'(\g^*)_+ =\{u|\, \l u,1\r=0\}$
is the space of distributions of integral $0$. Similarly, 
$\wedge\g^*=\R\oplus (\wedge\g^*)_+$ where $(\wedge\g^*)_+$
consists of elements of positive degree, and 
$$ \widehat{W}_G=\R\oplus (\widehat{W}_G)_+.$$
Let $\Pi:\,\widehat{W}_G\to \R\hra \wh{W}_G$ be projection 
defined by this splitting.
The differential $v^a \iota_a$ preserves the decomposition and
vanishes on $\R$, hence $\Pi$ is a chain map.

\begin{proposition}\label{Homotopy}
There exists an odd operator $h:\widehat{W}_G\to \widehat{W}_G$
satisfying
$$[h,v^a \iota_a]=\on{Id}-\Pi.$$
Hence $h$ provides a chain homotopy and 
$\widehat{W}_G$ is acyclic: $H(\widehat{W}_G)=\R$.
\end{proposition}

\begin{proof}
Under the pairing of $\wh{W}_G$ 
with $\Om(\g^*)$, the
projection $\Pi$ becomes dual to the projection map
$\pi:\,\Om(\g^*)\to \R\hra \Om(\g^*)$ induced by the inclusion of the
origin. Let $h^{Rh}:\,\Om^\star(\g^*)\to \Om^{\star-1}(\g^*)$ be the
standard de-Rham homotopy operator, so that 
\begin{equation} \label{RhamHomotopy}
[\d^{Rh},h^{Rh}]=\on{Id}-\pi. 
\end{equation}
Let $h:\,\wh{W}_G\to\wh{W}_G$ be minus the dual operator to $h^{Rh}$.
The Proposition follows by taking the dual of \eqref{RhamHomotopy}.
\end{proof}

An alternative proof of Proposition \ref{Homotopy}
can be given along the lines of Kumar-Vergne \cite{ku:eq},
Proposition 18.

\subsection{The non-commutative Weil algebra with generalized coefficients}
Just as for the usual Weil algebra it is important for applications
to introduce generalized coefficients in the non-commutative Weil algebra. 
For this we identify $u_a$ with the distribution on $G$, 
$$ u_a=\f{\d}{\d t}\Big|_{t=0} \delta_{\exp(t e_a)}, $$
and extend to an algebra homomorphism $U(\g)\to \ca{E}'(G)$. 
In this way $U(\g)$ is identified with the space of distributions 
on $G$ with support at $e$. 
The embedding $G\to \ca{E}'(G),\ g\mapsto \delta_g$ 
satisfies $\delta_{g_1} * \delta_{g_2}=\delta_{g_1g_2}$ where 
$u_1* u_2$ is the convolution of compactly supported distributions on $G$, 
defined as push-forward under group multiplication.  
For all $u\in \ca{E}'(G)$, $ \delta_g * u =(g\cdot u) * \delta_g$
which implies that the Lie derivative is given by a commutator 
$L_a u=[u_a,u]$.
We define the Weil algebra with generalized coefficients as
the super-algebra
$$ \widehat{\W}_G:=\ca{E}'(G)\otimes \Cl(\g).$$
The operators $\iota_a,L_a$ and $\d$ extend naturally to 
$\widehat{\W}_G$ and make it into a filtered $G$-differential algebra. 
As for the commutative Weil model, the multiplication map 
extends to a map
$$ \on{Mult}_\W:\widehat{\W}_{G\times G}\to \widehat{\W}_G$$
given as push-forward under group multiplication on 
the $\ca{E}'(G\times G)$ factor and by Clifford  
multiplication on the $\on{Cl}(\g\times\g)$ factor.
Given a $G$-differential algebra $B$ define the super-algebra 
$\widehat{\ca{H}}_G(B)$ as 
the cohomology of the complex $(\widehat{\W}_G\otimes B)_{basic}$,
$$ \widehat{\H}_G(B)= H((\widehat{\W}_G\otimes B)_{basic} ).$$
It is a module over the convolution algebra
$\widehat{\ca{H}}_G(\R)=\ca{E}'(G)^G$ 
of invariant distributions. As in the commutative case,
the embedding $U(\g) \hra \ca{E}'(G)$ induces
an embedding $\W_G \rightarrow \widehat{\W}_G$ and
a map $\ca{H}_G(B) \rightarrow \widehat{\H}_G(B)$ for
any $G$-differential algebra $B$.
For any compact oriented 
$G$-manifold $M$ the integration map descends to a map
$$\int_M:\,\widehat{\H}_G(\Om (M))\to \ca{E}'(G)^G.$$

The Cartan model with generalized coefficients is 
$$ \widehat{\ca{C}}_G(B)=(\ca{E}'(G)\otimes B)^G,$$
with differential and product structure given by the same formulas 
as in Section \ref{CartanModel}. Integration over $G$ defines a 
map
$$ \widehat{\ca{C}}_G(B)\to B^G$$
which becomes a chain map if $B^G$ is equipped with the new 
differential 
\begin{equation}
\label{TwistD}
\d+\f{1}{24} f_{abc} \iota_a \iota_b\iota_c.
\end{equation}
This map can also be described in the Weil model: Let 
$\Pi:\widehat{\ca{W}}_G\to \R$ the composition of 
the horizontal projection $P_{hor}:\, \widehat{\W}_G\to \ca{E}'(G)$
followed by integration over $G$. Then $\Pi$ defines a map
$$ \Pi\otimes 1:\,\widehat{\ca{W}}_G\otimes B\to B$$
which on basic elements agrees, under the isomorphism 
$P_{hor}\otimes 1:\ (\widehat{\ca{W}}_G\otimes B)^G\cong C_G(B)$
with the map defined above. In particular, the image 
of any basic element is closed under the differential 
\eqref{TwistD}.

\subsection{The pairing between $\wh{\W}_G$ and $\Om(G)$}
\label{Duality}
Identify $\Om^\star(G)=C^\infty(G)\otimes \wedge\g^*$ 
by means of left-invariant Maurer-Cartan forms $\theta_a$,
and identify $\wh{\W}_G\cong \ca{E}'(G)\otimes \wedge\g$ 
using the symbol map. In this section we study 
the pairing between $\wh{\W}_G$ and $\Om(G)$ given by these 
identifications.

We first describe the group analogue of the function
$\tau_0$. We will assume that the group $G$ is a direct product 
of a connected, simply connected group and a torus.

Let $\on{Spin}(\g)$ be the
spin group, defined as the image of 
$\on{so}(\g)\subset \on{Cl}^{(2)}_{even}(\g)$
under the (Clifford) exponential map. The $G$-action on $\Cl(\g)$ defines 
a homomorphism 
$$\g\to \on{so}(\g)\subset \Cl(\g),\,\mu\mapsto \mu_a g_a $$
where $g_a$ was defined in \eqref{Ga}. Using the assumption on $G$ 
it exponentiates to a map $ \tau:\,G\to \on{Spin}(\g)\hra \Cl(\g).$
Thus
\begin{equation}\label{DefTau}
\tau(\exp\,\mu)=\exp(-\hh f_{abc}\,\mu_a\,x_b\,x_c)
\end{equation}
which is similar to the definition of $\tau_0$ \eqref{DefTau0}.
By definition, $\tau(g_1)\tau(g_2)=\tau(g_1g_2)$ and $g\cdot x=
\Ad(\tau(g))x$ for all $x\in \Cl(\g)$. 
The two maps  $G\to \Cl(\g),\,g\mapsto\tau(g)$ and $G\mapsto \ca{E}'(G),\,
g\mapsto \delta_g$ combine into a map
$G\to \widehat{\W}_G,\ g\mapsto \tau(g)\delta_g$.
It has the property
\begin{equation}\label{ActionProperty}
\Ad(\tau(g)\delta_g)  w = g\cdot w 
\end{equation}
for all $g\in G$ and all $w\in \widehat{\W}_G$. 
\begin{proposition}
For all fixed $g\in G$ the element $\tau(g)\delta_g\in\widehat{\W}_G$ 
is closed.
\end{proposition}
\begin{proof}
Using the fact \eqref{ActionProperty} that conjugation by 
$\tau(g)\delta_g$ is the $G$-action on the Weil algebra, and 
that the element $\D$ is invariant, we find
$$ 
\d   \tau(g)\delta_g =
[\D, \tau(g)\delta_g]=
\tau(g)\delta_g\,(g^{-1}\cdot\ca{D}-\ca{D})
=0.
$$
\end{proof}
The function $\tau\in C^\infty(G)\otimes\Cl(\g)$ acts on $\W_G$ by
multiplication. This action commutes with $L_a$ since $\tau$ is
equivariant.  Let us conjugate the differential and contractions by
the action of $\tau$. Let $\theta_a$ and $\olt_a$ be the left/right
invariant Maurer-Cartan forms on $G$ and let
$$ \eta= \f{1}{12}f_{abc}\theta_a\theta_b\theta_c$$
be the canonical 3-form. It has the property
$\iota_a\eta=-\hh\d(\theta_a+\olt_a)$, from which one deduces 
that setting
$$\ti{\iota}_a=\iota_a+\f{1}{2}(\theta_a+\olt_a),\ \
\ti{L}_a=L_a,\ \ \ti{\d}=\d+\eta$$
gives $\Om(G)$ the structure of a $G$-differential space.
\begin{proposition}\label{Conj2}
The conjugates by $\tau$ of the Weil differential and contractions 
are dual to the differential $\ti{\d}$ and contractions 
$\ti{\iota}_a$ on $\Om(G)$:
$$
\Ad(\tau^{-1})(\d)=-(\d+\eta)^*,\ \
\Ad(\tau^{-1})(\iota_a)=-(\iota_a +\f{1}{2}(\theta_a+\olt_a))^*.
$$
\end{proposition}
\begin{proof}
The Weil differential is 
$$
\d=\ad(\D)=\D^L-\D^R=(x^L_a u_a^L-x_a^R u_a^R)+\ad(\gamma)
$$
Since $\gamma$ is $G$-invariant, $\Ad(\tau^{-1})\ad(\gamma)
=\ad(\gamma)$. 
We have 
$$\Ad(\tau^{-1})x_a^L=(\Ad_{g^{-1}})_{ab} x_b^L,\ \ 
\Ad(\tau^{-1})x_a^R=x_a^R
$$ 
and 
$$\Ad(\tau^{-1})u_a^L=u_a^L-(\Ad_{g^{-1}})_{ab} g_b^L,\ \ 
\Ad(\tau^{-1})u_a^R=u_a^R-g_a^L.$$
Hence 
\beq 
\lefteqn{\Ad(\tau^{-1})(u_a^L x_a^L-u_a^R x_a^R)}\\
&=&
(u_a^L-(\Ad_{g^{-1}})_{ab} g_b^L)(\Ad_{g^{-1}})_{ar} x_r^L
-(u_a^R-g_a^L)x_a^R\\
&=&\big(u_a^R x_a^L+\hh f_{abc}x_a^L x_b^L x_c^L\big)
-\big(u_a^R x_a^R+\hh f_{abc}x_a^L x_b^L x_c^R\big)\\ 
&=&u_a^R\iota_a+\hh f_{abc}x_a^L x_b^L \iota_a\\
&=&u_a^R\iota_a+\hh f_{abc}(y_a+\hh\iota_a) (y_b+\hh \iota_b) \iota_a.
\eeq
This result combines with the expression for $\ad(\gamma)$ 
from Proposition \ref{ClifDConj}
to 
$$ 
\Ad(\tau^{-1})(\d)=u_a^R \iota_a+\hh f_{abc}y_a\iota_b\iota_c
+\f{1}{12} f_{abc}\iota_a\iota_b\iota_c 
$$
which is indeed dual to $-(\d+\eta)^*$. The calculation for the 
contractions is simpler:
$$ \Ad(\tau^{-1})\iota_a=(\Ad_{g^{-1}})_{ab} {x}_b^L-x_a^R=
(\Ad_{g^{-1}}-I)_{ab}y_b+\hh(\Ad_{g^{-1}}+I)_{ab}\iota_b 
$$
which is dual to $-(\iota_a+\hh(\theta_a+\olt_a))$. 
\end{proof}

\begin{proposition}
The element 
$$ \Lambda=\exp(-x_a\olt_a)\tau(g)\delta_g
=\tau(g)
\exp(-x_a\theta_a)\delta_g
\in \wh{\W}_G\otimes\Om(G)
$$
has the properties, 
$$ \d\Lambda
=-\eta\Lambda,\ \
\iota_a\Lambda=-\hh(\theta_a+\olt_a)\Lambda.
$$
\end{proposition}

\begin{proof}
This follows from Proposition \ref{Conj2} since 
$\exp(-x_a\theta_a)\delta_g
\in \wh{\W}_G\otimes\Om(G)$ is the kernel of the identity map.
\end{proof} 

Let us also describe the co-multiplication on $\Om(G)$ which by
duality gives rise to the multiplication on $\wh{\W}_G$.

\begin{proposition}
Let $\on{Mult}_\W:\,\wh{\W}_G\otimes 
\wh{\W}_G\to\wh{\W}_G$ be the multiplication map for the Weil 
algebra. The composition  
$\tau^{-1}\circ \on{Mult}_\W\circ (\tau\otimes\tau):\,\wh{\W}_G\otimes 
\wh{\W}_G\to\wh{\W}_G$ is dual to the map 
$$ 
\exp(-\hh \theta_a^1\olt_a^2)\circ \Mult_G^*:\,\Om(G)\to \Om(G\times
G)$$
where $\on{Mult}_G:\,G\times G\to G$ is the group multiplication  
and the superscripts denote pull-backs to the respective $G$-factor.
\end{proposition}
\begin{proof}
Since $\tau^{-1}\La^1\La^2$ is the 
kernel of the multiplication map on $\wh{\W}_G$, 
it suffices to show that 
$$ 
\La^1 \La^2
=\exp({-\hh \theta_a^1\,\olt_a^2})(1\otimes\on{Mult}_G^*)\La 
.$$
Write $\La^1=\exp(-x_a \olt_a^1)\tau(g_1) \delta_{g_1}$ and similarly
for $\La^2$. We calculate: 
\beq
\La^1\La^2&=&\exp(-x_a \olt_a^1) \tau(g_1)\delta_{g_1}
             \exp(-x_b \olt_b^2)\tau(g_2)\delta_{g_2}\\
&=& \exp(-x_a \olt_a^1)
\exp(-x_b (\Ad_{g_1})_{cb} \olt^2_c)
\tau(g_1)\delta_{g_1}\tau(g_2)\delta_{g_2}\\
&=& \exp(-\hh \olt_a^1 (\Ad_{g_1})_{ba} \olt^2_b)\,
\exp(-x_a (\olt_a^1 +  (\Ad_{g_1})_{ba} \olt^2_b)
) \tau(g_1 g_2)\delta_{g_1 g_2}\\
&=& \exp(-\hh \theta_a^1 \olt_a^2) (1\otimes \on{Mult}_G^*)\La
\eeq
For the last equality we used that 
$$ \on{Mult}_G^*\olt_a = \olt_a^1 + (\Ad_{g_1})_{ba}\olt_b^2,$$
and for the third equality we used that 
$$ \exp(\kappa_a^1 x_a)\exp(\kappa_a^2 x_a)=
\exp(-\hh \kappa_a^1\kappa_a^2)
\exp((\kappa_a^1+\kappa_a^2) x_a)
$$
if $\kappa_a^1,\kappa_a^2$ are elements of a commutative 
super-algebra.
\end{proof}

\section{The quantization map $\ca{Q}:\,W_G\to \W_G$}
\label{D}
In this Section we construct an explicit isomorphism of 
$W_G$ and $\W_G$ as $G$-differential spaces. Notice that by 
contrast, the exterior algebra $\wedge\g^*$ and the Clifford
algebra $\Cl(\g)$ are not isomorphic as differential spaces
(unless $\g$ is abelian): As we have seen the 
cohomology of $(\Cl(\g),\ad(\gamma))$ is trivial in all degrees 
while the cohomology $(\wedge\g,\d)$ for the Lie algebra differential
is not.

\subsection{The Duflo map}
The Birkhoff-Witt symmetrization map $S(\g)\to U(\g)$ is 
the unique $G$-module isomorphism sending $(\mu_a v_a)^k$ to 
$(\mu_a u_a)^k$, for all $k\ge 0$ and all $\mu\in\g$.
Under the identification of $U(\g)$ with distributions on $G$
supported at the identity, and $S(\g)$ with distributions on $\g$
supported at $0$, this isomorphism is induced by 
push-forward under the exponential map:
$\exp_*:\, \ca{E}'(\g)\to\ca{E}'(G)$.  

Duflo \cite{du:op} introduced
a different isomorphism $S(\g)\cong U(\g)$ given by composition
of $\exp_*$ with multiplication by the square root of the Jacobian of
the exponential map, $J^\hh:\,\g\to\R$. (Recall that the Jacobian $J$
has a globally defined smooth square root $J^\hh$ with $J^\hh(0)=1$.)
The Duflo map
$$ \on{Duf}:=\exp_*\circ J^\hh $$
has the important property that it induces an {\em algebra isomorphism} 
$S(\g)^G\to U(\g)^G$. 
We will also refer to the map $\on{Duf}:=\exp_*\circ
J^\hh:\,\ca{E}'(\g)\to \ca{E}'(G)$ as the Duflo map; it restricts to a
ring homomorphism for the convolution algebras, $\ca{E}'(\g)^G\to
\ca{E}'(G)^G$. 

Combining the Duflo map with the inverse of the symbol 
map, and using the inner product to identify $\g\cong\g^*$, 
we obtain a map between Weil algebras
$$(\on{Duf}\times\sig^{-1}):\,\widehat{W}_G \to \widehat{\W}_G.$$
In the following Section \ref{WeilQuantization} 
we show that one can do much better: 
There exists a map $\ca{Q}: \wh{W}_G \rightarrow \wh{\W}_G$ 
which we call the {\em quantization map} such that
$\ca{Q}$ is a homomorphism of $G$-differential spaces.

\subsection{Quantization of Weil algebra}\label{WeilQuantization}
The definition of the quantization map $\ca{Q}: \wh{W}_G \rightarrow 
\wh{\W}_G$
involves a certain skew-symmetric 
tensor field $T_{ab}$, given as follows. 
Let $\g_\#\subset\g$ be the open subset on which the exponential map is 
a local diffeomorphism, that is $\g_\#=\g\backslash J^{-1}(0)$.
Let $e_a^L$ be the left-invariant vector field on $G$ which 
equals $e_a$ at the group unit, and $e_a^R$ the corresponding 
right-invariant vector field.
It turns out (see Appendix A) that at any point $\mu\in\g_\#$, 
the pull-back of the half-sum $\f{e_a^L+e_a^R}{2}$ under the exponential 
map differs from the constant vector field $\f{\p}{\p\mu_a}$ 
only by a vector tangent to the $G$-orbit through $\mu$. 
It follows that there is a unique tensor field 
$T:\,\g_\#\to\g\otimes\g$  such that
$T(\mu)$ takes values in $\g_\mu^\perp\otimes \g_\mu^\perp$ and 
\begin{equation}\label{Half}
\exp^*\Big(\f{e_a^L+e_a^R}{2}\Big)-\f{\p}{\p\mu_a}=T_{ab} (e_b)_\g.
\end{equation}
We verify in Appendix A that $T$ is explicitly given as 
$$ T\in C^\infty(\g_\#,\g\otimes\g),\
T_{ab}(\mu):=f(\ad_\mu)_{ab}.$$
where $f(s)$ is the function  
$$ f(s)=\f{1}{s}-\hh\on{cotanh}\big(\f{s}{2}\big).$$
In particular $T_{ab}=-T_{ba}$.  
The product $J^\hh \exp\Big(-\hh T_{ab} \iota_a\iota_b\Big)$ is a 
smooth function on $\g$ with values in operators on 
$\wh{W}_G$, since the zeroes of $J^\hh$
compensate the singularities of $T$.
The main result of this paper is presented in the following Theorem.
\begin{theorem}\label{Duflo}
The {\em quantization map}
$$\ca{Q}:=(\on{Duf}\times\sig^{-1})\circ \exp
\Big(-\hh T_{ab}\iota_a\iota_b\Big):\,\widehat{W}_G 
\to \widehat{\W}_G$$
is a homomorphism of $G$-differential spaces. That is, $\ca{Q}$ 
satisfies 
$$\ca{Q}\circ L_a=L_a\circ \ca{Q},\ \  
\ca{Q}\circ \iota_a=\iota_a\circ \ca{Q},\ \  
\ca{Q}\circ \d=\d\circ \ca{Q}.$$
\end{theorem}

\begin{proof}
The properties $\ca{Q}\circ L_a=L_a\circ \ca{Q}$
and $\ca{Q}\circ \iota_a=\iota_a\circ \ca{Q}$ is obvious.
Let us show that $\ca{Q}$ is a chain map. For clarity 
we denote the differential on $\wh{\W}_G$ by $\d^\W$ 
and the differential on $\wh{W}_G$ by $\d^W$. 

We want to compare
the expression for $\d^\W$ obtained in Proposition
\ref{ConjSigma} to
\begin{equation}\label{ConjugateD}
\Ad\Big(J^\hh \exp(-\hh T_{ab}\iota_a\iota_b)\Big)\d^W.  
\end{equation}
Using the invariance property
$ [L_j\otimes 1,T_{ab}]+f_{jar} T_{rb}+ f_{jbs} T_{as}=0$
we compute:
\beq
\lefteqn
{\ad(-\hh T_{ab}\iota_a\iota_b)\d^W}\\
&=&
-\hh \f{\p T_{rs}}{\p \mu_a} \iota_r\iota_s\iota_a
+\f{1}{4} T_{rs} f_{abc}[\iota_r\iota_s,y_b y_c] \iota_a
-T_{rs} \iota_r [\iota_s,y_j]\,L_j\otimes 1 -
\hh [T_{rs},L_j\otimes 1]\,y_j\iota_r\iota_s
\\
&=& - \hh  \f{\p T_{rs}}{\p \mu_a} \iota_r\iota_s\iota_a -
f_{ajk} T_{rj} y_k \iota_r \iota_a
-\hh T_{bc} f_{abc}\iota_a
+f_{jrk}T_{ks}\,y_j\iota_r\iota_s-
T_{rj} \iota_r (L_j\otimes 1)\\
&=&-\hh T_{bc} f_{abc}\iota_a - \hh  \f{\p T_{rs}}{\p \mu_a} \iota_r\iota_s\iota_a
-T_{rj} \iota_r (L_j\otimes 1),\\
\lefteqn{\f{1}{2!}\ad^2(-\hh T_{ab}\iota_a\iota_b)\d^W}\\ &=&
\f{1}{4} T_{rj}\iota_r\iota_u\iota_v [T_{uv}, L_j\otimes 1]\\
&=&-\hh T_{rj} T_{kv} f_{juk} \iota_r\iota_u\iota_v.
\eeq
It follows that  
$\Ad(\exp(-\hh T_{ab}\iota_a\iota_b))\d^W$ is given by the formula
$$
\d^W -T_{ab} (L_b\otimes 1) \iota_a  -\hh\Big(
\f{\p T_{ab}}{\p \mu_c}+T_{ar} f_{rbs} T_{sc} \Big) \iota_a\iota_b \iota_c
-\hh f_{abc}T_{bc} \iota_a.
$$
In Appendix A we show that the tensor field $T_{ab}$ is a solution of the
classical dynamical Yang-Baxter equation 
$$ \on{Cycl}_{abc}\Big(\f{\p T_{ab}}{\p \mu_c} +
T_{ar} f_{rbs} T_{sc}\Big)=\f{1}{4}f_{abc}$$
where $\on{Cycl}_{abc}$ means the sum over cyclic permutations of the 
indices $a,b,c$. Therefore, 
$$ \Ad\big(e^{-\hh T_{ab}\iota_a\iota_b}\big)(\d^W)=
\d^W -T_{ab} (L_b\otimes 1) \iota_a  
-\f{1}{24}f_{abc} \iota_a\iota_b \iota_c
-\hh f_{abc} T_{bc} \iota_a.
$$
Conjugating further by $J^\hh$ adds $\hh\f{\p \ln J}{\p \mu_a}$ 
to this expression. As shown in Appendix \ref{AppA}, the 
derivatives of $J$ can be expressed in terms of the 
tensor field $T_{ab}$ as follows:
$$ f_{abc}T_{bc}-\f{\p \ln J}{\p \mu_a}=0.$$
Hence, 
\beq
\lefteqn{\Ad\Big(J^\hh e^{-\hh T_{ab}\iota_a\iota_b}\Big)\d^W}\\
&=&
  \d^W -T_{ab} (L_b\otimes 1) \iota_a
-\f{1}{24} f_{abc} \iota_a\iota_b \iota_c\\
&=&y_a(L_a\otimes 1)+
\Big(v_a-T_{ab} (L_b\otimes 1)-\hh f_{abc} y_b y_c\Big) \iota_a  
-\f{1}{24} f_{abc}\iota_a\iota_b \iota_c.
\eeq
The definition of $T$ implies that  
$$ \f{u_a^L +u_a^R}{2}\circ \exp_* = \exp_*\circ 
(v_a - T_{ab} (L_b\otimes 1))$$
Together with Proposition \ref{ConjSigma} this shows that 
$$ \exp_*\circ \Ad\Big(J^\hh e^{-\hh T_{ab}\iota_a\iota_b}\Big)\d^W
=\d^\W\circ \exp_*,$$
which means that $\ca{Q}$ is a chain map.
An alternative proof of this fact will be given 
in Section \ref{transpose}.
\end{proof}

\begin{remark}
Notice that the quantization map restricts to the 
Duflo map $\on{Duf}:\,\ca{E}'(\g)\to \ca{E}'(G)$ on the subalgebra 
$\ca{E}'(\g)\otimes 1\subseteq \widehat{W}_G$,   
and to the inverse of the symbol map $\sig$ on $1\otimes \wedge\g$.
Furthermore, $\ca{Q}$ restricts to an isomorphism of $G$-differential 
spaces $\ca{Q}:\,W_G\to \W_G$.
\end{remark}

As a direct consequence to Theorem \ref{Duflo}, we have:
\begin{corollary}
For any $G$-differential algebra $B$, the quantization map 
$ \ca{Q} $ induces a linear isomorphism ${H}_G(B)\cong {\H}_G(B)$
and a linear map
$\widehat{H}_G(B) \rightarrow  \widehat{\H}_G(B)$. 
\end{corollary}

\subsection{Quantization map in Cartan model}
Suppose that $B$ is a $G$-differential algebra. The following 
proposition describes the  chain map 
$$(\widehat{W}_G\otimes B)_{basic}\to (\widehat{\W}_G\otimes B)_{basic}.$$
induced by $\ca{Q}$ in terms of the Cartan model:
\begin{proposition}
The chain map between the Cartan models 
$$ (\ca{E}'(\g)\otimes B)^G \to (\ca{E}'(G)\otimes B)^G$$
induced by $\ca{Q}$ is given by 
$$ \on{Duf}\circ\exp\Big(-\hh T_{ab} (1\otimes \iota_a\iota_b)\Big).$$ 
\end{proposition} 
\begin{proof}
We have to show that on the subspace of 
basic elements $(\widehat{W}_G\otimes B)_{basic}$, 
$$ (P_{hor}\otimes 1)\circ \ca{Q}=
\on{Duf}\circ\exp
\Big(-\hh T_{ab} (1\otimes \iota_a\iota_b)\Big)\circ (P_{hor}\otimes 1)$$
Since basic elements are in particular horizontal, the operator 
$\exp
\big(-\hh T_{ab} (\iota_a\iota_b\otimes 1)\big)$ appearing in 
the definition of $\ca{Q}$ can be replaced with 
$\exp
\big(-\hh T_{ab} (1\otimes \iota_a\iota_b)\big)$, which then commutes with 
$P_{hor}\otimes 1$. 
\end{proof}

\section{The transpose of the quantization map}
\label{transpose}

Let $\varpi\in\Om^2(\g)$ be the image of the closed form 
$\exp^*\eta\in\Om^3(\g)$ under the de Rham homotopy 
operator $\Om^\star(\g)\to \Om^{\star-1}(\g)$.
Explicitly,
$$\varpi_\mu=-\hh g(\ad_\mu)_{ab}\d\mu_a\d\mu_b$$ 
where $g(s)$ is the function 
$$ g(s)=\f{\sinh(s)-s}{s^2}.$$

\begin{theorem}\label{duality}
Under the pairings of $\wh{W}_G$ with $\Om(\g)$ and of 
$\wh{\W}_G$ with $\Om(G)$, the composition 
$\tau^{-1}\circ \ca{Q}\circ \tau_0$ is dual to the map 
$$ e^\varpi\circ \exp^*:\,\Om(G)\to \Om(\g)$$
\end{theorem}

Since this map is a chain map for the differentials 
$\d$ on $\Om(\g^*)$ and $\d+\eta$ on $\Om(G)$, this 
gives an alternative proof of the fact that $\ca{Q}$ 
is a chain map.
To prove this result note that $\tau^{-1}\ca{Q}(\Lambda_0)\in \wh{\W}_G
\otimes\Om(\g^*)$ is the integral kernel of the quantization map,
while $\tau^{-1}\exp^*\Lambda$ is the kernel of the map 
$\exp_*:\,\wh{W}_G\to \wh{\W}_G$. 
Hence, the theorem will follow once we 
show that $\ca{Q}(\Lambda_0)=e^\varpi\,\exp^*\Lambda$. 

We begin by computing the quantization of the form 
$\tau_0(\mu)\delta_\mu$. This will require the 
following Lemma. 

\begin{lemma}\label{BerezinLemma}
Let $V$ be an oriented Euclidean vector space of even dimension 
$\dim V=2n$, 
and suppose $S\in \on{so}(V)$ is invertible. Let $e_1\ldots e_{2n}$ 
be an oriented orthonormal basis of $V$.
Then 
$$ \f{1}{\Pf(S)}\exp(\hh (S^{-1})_{ab}\iota_a \iota_b) 
\exp(\hh S_{ab} e_a \wedge e_b) 
=e_1\wedge\ldots\wedge e_{2n}$$
where $\on{Pf}(S)$ is the Pfaffian of $S$.
\end{lemma}

\begin{proof}
Block-diagonalizing $S$ one reduces the Lemma to the case $\dim V=2$
and $S_{11}=S_{22}=0,\ S_{12}=-S_{21}=s$.
In this case 
the equation becomes 
$$\f{1}{s}\exp(\f{1}{s}\iota_1\iota_2)\exp(s e_1\wedge e_2)= e_1\wedge e_2$$
which follows immediately by writing $ \exp(s e_1\wedge e_2) =1+s
e_1\wedge e_2$ and
$\exp(\f{1}{s}\iota_1\iota_2)=1+\f{1}{s}\iota_1\iota_2$.
\end{proof}

\begin{proposition}\label{Qtau}
For all $\mu\in\g$, 
$$ \ca{Q}\,\, \tau_0(\mu)\delta_\mu   =\tau(\exp\mu)\delta_{\exp\mu}  .$$
\end{proposition}

\begin{proof}
Recall that the set of regular elements $G_{reg}$ is the set
of all elements whose stabilizer is a maximal torus.
Since both sides depend continuously on $\mu$ we may assume 
$\mu\in\exp^{-1}(G_{reg})$.
  
Choose any orientation of $\g_\mu^\perp$. We will apply the Lemma with 
$V=\g_\mu^\perp$.
Let $\Pf_{\g_\mu^\perp}$ denote the Pfaffian of an operator 
on the orthogonal complement of $\g_\mu=\ker(\ad_\mu)$
(using the metric induced from $\g$ and  any choice of orientation).
The operator $\ad_\mu$ has kernel $\ker(\ad_\mu)={\g_\mu}$ and is 
invertible on ${\g_\mu}^\perp$. The 
square root of the Jacobian of the exponential map is a 
quotient of two Pfaffians (see e.g. \cite{he:di}, p. 105):
$$ J^\hh(\mu)=\f{\det^\hh_{{\g_\mu}^\perp}(2\sinh(\f{\ad_\mu}{2})}
{\det^\hh_{{\g_\mu}^\perp}({\ad_\mu})}.
$$
Consider the skew-symmetric tensor fields
$\mathfrak{r}_0\in C^\infty(\g_{reg},\g\wedge\g)$ and
$\mathfrak{r}\in C^\infty(G_{reg},\g\wedge\g)$ on the set of regular
elements $\g_{reg}$ resp. $G_{reg}$
introduced in the Appendix. 

We can re-write the quantization map as 
$$ \ca{Q}=(\exp_*\times  \sig^{-1})\ 
{{\det}^\hh(\cosh(\f{\ad_\mu}{2}))}
\ca{T}^{-1}\ \ca{T}_0 $$ 
where 
$$ \ca{T}_0=\f{1}{\det^\hh_{{\g_\mu}^\perp}(\f{\ad_\mu}{2})}
\exp(-\hh (\mathfrak{r}_0)_{ab} \iota_a\iota_b)
$$ 
and 
$$ \ca{T}=\f{1}{\det^\hh_{{\g_\mu}^\perp}(\tanh(\f{\ad_\mu}{2}))}
\exp(-\hh \mathfrak{r}_{ab} \iota_a\iota_b).
$$ 
As a special case of \cite{be:he}, Proposition 3.13, 
the symbol of $\tau$ is given by 
\begin{equation}\label{symboltau} 
\sig(\tau(g))=
{\det}^\hh(\cosh(\f{\ad_\mu}{2}))
\exp\Big(-\,\big(\f{\Ad_g-1}{\Ad_g+1}\big)_{ab}
y_a y_b\Big)
\end{equation}
where $g=\exp(\mu)$.
Together with Lemma \ref{BerezinLemma} 
this shows 
$ {\det}^\hh(\cosh(\f{\ad_\mu}{2}))\ca{T}_0\tau_0= \ca{T}\sig(\tau)$.
\end{proof}

\begin{proposition}\label{ImageQuantization}
The image of $\La_0$ under the quantization 
map is given by 
$$ (\ca{Q}\otimes 1)\La_0
=\exp(\varpi)\ (1\otimes\exp^*)\La. $$
\end{proposition}

\begin{proof}
By continuity it suffices to verify the equation on the open dense 
subset $\g_{reg}$ of regular elements of $\g$.  
We define 1-forms $\kappa_a\in \Om^1(\g_{reg})$ by 
$$ \kappa_a := (\mathfrak{r}_0)_{ab} \d\mu_b. $$
At any point $\mu\in\g_{reg}$ we have a decomposition of the tangent
space into the spherical part, i.e. the tangent space to the orbit
$(T_\mu\g)^{sp}=\on{im}(\ad_\mu)$, and its orthogonal complement, the
radial part $(T_\mu\g)^{rad}=\on{ker}(\ad_\mu)$ (spanned 
by the Cartan subalgebra containing $\mu$). Correspondingly 
every vector field on $\g_{reg}$, and dually every 1-form, decomposes 
into radial and spherical part. The spherical part is spanned by the forms 
$\kappa_a$. We have 
$$ (\d\mu_a)^{sp}=(\ad_\mu)_{ab}\kappa_b,\ \  (\exp^*\olt_a)^{sp}=(1-e^{-\ad_\mu})_{ab}\kappa_b
$$
while 
$$ (\d\mu_a)^{rad}=(\exp^*\olt_a)^{rad}.$$
Then the definition of $\varpi$ can be re-written
$\varpi=\varpi_2-\varpi_1$ where 
\beq \varpi_2 &=&-\hh (\sinh(\ad_\mu))_{ab}\kappa_a \kappa_b, \\
\varpi_1&=&  -\hh (\ad_\mu)_{ab} \kappa_a \kappa_b.
\eeq
The form $\varpi_1\in\Om^2(\g_{reg})$ combines nicely with 
$\La_0$: 
\begin{eqnarray}
\lefteqn{\exp(\varpi_1)\La_0=
\exp(\varpi_1-y_a d\mu_a- \hh f_{abc} \mu_a y_b y_c)\,\delta_\mu}
\nonumber\\
&=&\exp(-y_a (d\mu_a)^{rad}  ) 
\exp(-\hh f_{abc} \mu_a \,(y_b+\kappa_b) (y_c+\kappa_c))\delta_\mu
\nonumber \\
&=&\exp(-y_a (d\mu_a)^{rad}) 
\exp(\kappa_a\iota_a) \tau_0(\mu)\delta_\mu.\label{One}
\end{eqnarray}
We will see that there exists a similar expression for 
$\exp(\varpi_2)\exp^*\La$. We will need the following 
\begin{lemma}\label{DiLemma}
Let $V$ be an oriented Euclidean vector space, and $S\in \on{so}(V)$. 
Let $x_a\in \on{Cl}(V)$ be the generators corresponding to some choice
of oriented, orthonormal basis of $V$, 
and let $\kappa_a$ be odd elements in some 
commutative super-algebra $\ca{A}$. Then the following identity 
in $\Cl(V)\otimes\ca{A}$ holds: 
\begin{equation}\label{Cl1}
\exp(-\iota_a \kappa_a)
\exp\big(\hh S_{ab}x_a x_b\big)
= 
\exp(\varpi_2)  \exp(-x_r\gamma_r) 
\exp\big(\hh S_{ab}x_a x_b\big) 
\end{equation}
where $\gamma_r=(1-e^{S})_{rs}\kappa_s$ and 
$\varpi_2=\hh (\on{sinh}(S))_{ab}\kappa_a\kappa_b$. 
\end{lemma}
\begin{proof}
Block-diagonalizing $S$ we can assume 
$\dim V=2$ and that $S_{12}=-S_{21}=s$ and
$S_{11}=S_{22}=0$. Then 
$$ \exp(\hh S_{ab}x_a x_b)=
\exp(s x_1 x_2)=
\cos(\f{s}{2})+2\sin(\f{s}{2})x_1 x_2.$$
Furthermore, $\varpi_2=\sin(s)\kappa_1\kappa_2$ and 
\beq
\gamma_1&=&
2\sin(s/2)(\sin(s/2)\kappa_1-\cos(s/2)\kappa_2),\\
\gamma_2&=&
2\sin(s/2)(\sin(s/2)\kappa_2+\cos(s/2)\kappa_1).
\eeq
Equation \eqref{Cl1} becomes
$$
(1+\kappa_1 \iota_1 +\kappa_2 \iota_2
-\kappa_1\kappa_2\iota_1\iota_2)(\cos(s/2)+2\sin(s/2)x_1 x_2)
$$
$$=
(1+\sin(s)\kappa_1\kappa_2)
(1-x_1 \gamma_1 - x_2 \gamma_2- x_1 x_2\gamma_1\gamma_2)
(\cos(s/2)+2\sin(s/2)x_1 x_2)
$$
which is verified by an elementary calculation.
\end{proof}
Using the Lemma, we can write
\begin{equation}\label{Two}
\exp(\varpi_2) \exp^*\Lambda=
\exp(-x_a \exp^*
\olt_a^{rad})\exp(\kappa_a\iota_a)\tau(\exp\mu)\delta_{\exp\mu}.
\end{equation}
Since the quantization map is equivariant, 
and since it intertwines the 
contractions on $\widehat{W}_G$ and $\widehat{\W}_G$, 
Equations \eqref{One}, \eqref{Two} show that 
$$ \ca{Q}(e^{\varpi_1}\Lambda_0)=e^{\varpi_2}\exp^*\Lambda.$$
\end{proof}

\section{Ring structure}\label{Kernels}
In the Section \ref{D} we established that the quantization 
map induces 
maps in cohomology $\ca{Q}:\,\wh{H}_G(B)\to \wh{\H}_G(B)$
and $\ca{Q}:\,{H}_G(B)\to {\H}_G(B)$. 
Are these maps algebra homomorphisms? 

To answer this question let us introduce two 
homomorphisms of $G$-differential algebras.
$\phi_1,\phi_2:\,\wh{W}_G\otimes\wh{W}_G\to \wh{\W}_G$ 
by 
\beq \phi_1&=&\on{Mult}_\W\circ (\ca{Q}\otimes \ca{Q}),\\
     \phi_2&=&\ca{Q}\circ \on{Mult}_W.
\eeq
\begin{theorem}\label{Duflo1}
The maps $\phi_1,\phi_2:\,\wh{W}_G\otimes\wh{W}_G\to \wh{\W}_G$ are
$G$-chain homotopic (cf. Definition \ref{ChainHomotopic}). They restrict to
$G$-chain homotopic maps $W_G\otimes W_G\to \W_G$. It follows 
that 
for any $G$-differential algebra $B$, the quantization map 
induces a ring isomorphism $H_G(B)\to \ca{H}_G(B)$ and a ring 
homomorphism $\wh{H}_G(B)\to \wh{\H}_G(B)$. 
\end{theorem}
Taking $B$ to be the trivial 
$G$-differential algebra $B=\R$ we recover the fact that 
the Duflo map induces a ring homomorphism 
$\ca{E}'(\g)^G\to \ca{E}'(G)^G$ and a ring isomorphism
$(S\g)^G\cong U(\g)^G$. 

\begin{proof}
For $j=1,2$ let 
$\ca{K}_j\in \wh{\W}_G\otimes\Om(\g\times\g)$
be defined as 
$$ \ca{K}_j=\phi_j(\Lambda_0\otimes\Lambda_0).$$
then $\tau^{-1}\ca{K}_j$ are integral kernels for 
the maps $\phi_j$. From the fact that these are chain maps
intertwining contractions, one obtains
$$ \d\ca{K}_j=0,\ \ \iota_a\ca{K}_j=-\on{Add}^*_\g(\d\mu_a)\ca{K}_j.
$$ 
Since $\on{Mult}_W(\Lambda_0\otimes\Lambda_0)$
is just $(1\otimes \on{Add}_\g^*)\Lambda_0$ (pull-back under the addition map
$\on{Add}_\g:\,\g\times\g\to \g$), 
$$ \ca{K}_2=\on{Add}_\g^*\ca{Q}(\Lambda_0)=
\on{Add}_\g^*(e^\varpi\exp^*\Lambda).$$
Notice that $e^\varpi\exp^*\Lambda$ is an {\em invertible} element 
of the algebra $\wh{\W}_G\otimes \Om(\g)$, and the same is true 
for its pull-back under $\on{Add}_\g$. Hence $\ca{K}_2^{-1}$ is 
defined, and the product 
$$ \ca{N}=\ca{K}_2^{-1}\ca{K}_1\in  \wh{\W}_G\otimes \Om(\g\times\g)$$
is $G$-basic.
Since $\phi_j(1\otimes 1)=1$ for both $j=1,2$, the 
pull-back of both $\ca{K}_j$, hence also of 
$\ca{N}$ to the origin in $\g\times\g$ is the 
identity element $1\in\wh{W}_G$. 
Since $\g\times\g$ retracts equivariantly to the origin, 
it follows that there exists an odd element 
$\Gamma\in (\wh{\W}_G\otimes \Om(\g\times\g))_{basic}$ 
with $\ca{N}=1+\d\Gamma$. Multiplying this identity by $\ca{K}_2$
we find
$$ \ca{K}_2=\ca{K}_1+\d(\ca{K}_2\Gamma).$$
The element $\ca{K}_2\Gamma$ has the property 
$$ \iota_a(\ca{K}_2\Gamma)=-\d\mu_a\,(\ca{K}_2\Gamma).$$
Consequently, $\tau^{-1}(\ca{K}_2\Gamma)$ is the kernel
for an odd linear map
$$h:\,\wh{W}_G\otimes\wh{W}_G\to \wh{\W}_G $$
that provides a $G$-chain homotopy between $\phi_1$ and $\phi_2$. 
We notice that $h$  restricts to a $G$-chain homotopy  
$W_G\otimes W_G\to \W_G$ between the restrictions of 
$\phi_1$ and $\phi_2$.
\end{proof}

\begin{remark}
The first part of Theorem \ref{Duflo1}
shows more generally that if 
$B$ is a $G\times G$-differential space, the natural diagram
\begin{equation*}\
\vcenter{\xymatrix{  
{ \widehat{H}_{G\times G}(B)}
\ar[d] 
\ar[r]& {\widehat{H}_G(B)}\ar[d]\\
{\widehat{\ca{H}}_{G\times G}(B) } 
\ar[r]
& { \widehat{\ca{H}}_G(B) }
}}
\end{equation*}
where the vertical arrows are quantization maps and the 
horizontal arrows are induced by multiplication in the 
Weil algebras, commutes. 
\end{remark}

The simplest non-trivial example illustrating the ring isomorphism 
$H_G(B)\cong \H_G(B)$
is the rotation action of $G=\SU(2)$ on the 2-sphere
$S^2\subset \R^3$. As is well-known, the equivariant cohomology 
$H_{G}(\Om(S^2))$ as 
an algebra over $\R$ is a polynomial ring in one generator
$[\omega]$ of degree 2. We can normalize $[\om]$ by 
the condition $[\om]^2=[\lambda]$ where 
$\lambda=v_a v_a \in S(\g)$ 
is the generator of the ring of invariant polynomials.
Let us verify explicitly that $\ca{Q}([\om]^2)=\ca{Q}([\om])^2$.

To describe a representative for $[\om]$ let us choose the basis 
of $\g$ according to the standard identification 
$\mathfrak{su}(2)\cong\R^3$, so that the structure constants are given 
by the totally anti-symmetric tensor $f_{abc}=\varepsilon_{abc}$.
Also, let us view $S^2$ as the unit sphere in $\R^3$, with coordinates
$n_a$. In the Cartan model, a representative for $[\om]$ is given by
$$ \omega= \frac{1}{2} f_{abc} n^a dn^b dn^c +
n^a v^a.$$ 
The representatives satisfy 
$$ \omega^2 = \lambda + d_G(f_{abc} v^a n^b d n^c). $$
The image of $\lambda$ under the quantization map, 
which on $(S\g)^G$ restricts to the Duflo map $\on{Duf}$,
is the quadratic Casimir, 
$$ \ca{Q}(\lambda)=u^a u^a + \frac{1}{4} \in U(\g)^G. $$
Applying the quantization map to $\omega$, one obtains
$$ \ca{Q}(\omega) =\frac{1}{2} f_{abc} n^a dn^b dn^c +
n^a u^a. 
$$
Using 
$$
\frac{1}{8} \sum_{\alpha,\beta}(\iota_\alpha \iota_\beta \ca{Q}(\omega) )^2 
=\frac{1}{4} n^a n^a =
\frac{1}{4}, 
$$
we find that the square of $\ca{Q}(\omega)$ with respect to the ring structure 
$\odot$ is 
\beq \ca{Q}(\omega) \odot \ca{Q}(\omega)&=& \ca{Q}(\omega)^2 + 
\frac{1}{2}\sum_\alpha (\iota_\alpha \ca{Q}(\omega))^2 + 
\frac{1}{8} \sum_{\alpha,\beta}(\iota_\alpha \iota_\beta \ca{Q}(\omega))^2 \\
&=&
u^au^a + d_G(f_{abc} u^a n^b dn^c) 
+\frac{1}{4}\\
&=&\ca{Q}(\lambda)+d_G(f_{abc} u^a n^b dn^c)
=\ca{Q}(\om^2).
\eeq

\begin{appendix}
\section{Properties of the tensor $T_{ab}$}\label{AppA}
In this Section we prove the properties of the tensor field
$T$ which were used in this paper. Let us first verify that the  
definition $T(\mu)=f(\ad_\mu)$, where 
$$ f(s)=\f{1}{s}-\hh\on{cotanh}\big(\f{s}{2}\big),$$
coincides with the definition given in the text.
\begin{lemma}\label{Lem2}
The half-sum of the left- and the right-invariant vector fields 
satisfies 
\begin{equation}\label{Halfsum}
\exp^*\  \f{e_a^L +e_a^R}{2}=\f{\p}{\p \mu_a} + T_{ab} (e_b)_\g.
\end{equation}
\end{lemma}
\begin{proof}
Consider the function $g(s)=\f{1-e^{-s}}{s}$. The right-invariant 
Maurer-Cartan form $\olt$ has the property 
$\exp^*\olt_a=(g(\ad_\mu))_{ab}\d\mu_b$. Hence the right-invariant 
vector field satisfies 
$$\exp^*(e_a^R)=
(g(\ad_\mu))_{ba}\f{\p}{\p \mu_b}=
(\check{g}^{-1}(\ad_\mu))_{ab}\f{\p}{\p \mu_b}.$$
where $\check{g}(s)=g(-s)$. Similarly 
$$ \exp^*(e_a^L)=(g^{-1}(\ad_\mu))_{ab}\f{\p}{\p \mu_b}.$$
\, From this the Lemma follows since $1+s f=(g^{-1}+\check{g}^{-1})/{2}.$
\end{proof}
\begin{lemma}\label{Lem1}\label{DerivativeJ}
The derivative of the Jacobian of the exponential map is given 
by the formula
$$  \f{\p \ln J}{\p \mu_a} =f_{abc} T_{bc}.$$ 
\end{lemma}
\begin{proof}
It is well-known (see e.g. \cite{he:di}, p.105) 
that $J(\mu)=\det(g(\ad_\mu))$
with $g$ as in the proof of the previous Lemma.
Therefore
\beq 
\f{\p \ln J}{\p \mu_a}&=&\f{\p }{\p \mu_a}\ln (\det(g(\ad_\mu)))\\
&=&\f{\p }{\p \mu_a} \on{tr}(\ln(g(\ad_\mu))) \\
&=&\on{tr} \Big((\ln g)'(\ad_\mu) \f{\p }{\p \mu_a}(\ad_\mu)\Big)\\
&=& f_{abc}\Big((\ln g)'(\ad_\mu)\Big)_{bc}
\eeq
Since the anti-symmetric part of $(\ln g)'(s)=(e^{s}-1)^{-1}-{s}^{-1}$ is 
equal to $f$ the final result is $f_{abc} T_{bc}$. 
\end{proof}

The most interesting and most 
complicated property of $T$ is the following result:
\begin{lemma}[Etingof-Varchenko \cite{et:ge}] 
\label{Lem3}  
The tensor $T$ is a solution of the classical dynamical Yang-Baxter 
equation 
with coupling constant
$1/4$:
\begin{equation}\label{YangBaxter}
 \on{Cycl}_{abc}(\f{\p T_{bc}}{\p \mu_a}+T_{ar}f_{rbs}T_{sc})
=\f{1}{4}f_{abc},
\end{equation}
where $\on{Cycl}_{abc}$ denotes the sum over cyclic permutations of 
$a,b,c$. 
\end{lemma}
This solution of  the classical dynamical Yang-Baxter equation was 
obtained by Etingof-Varchenko as part of their general classification 
scheme, see \cite{et:ge}, Theorems 3.1 and 3.14. 
The remainder of this section is devoted to a 
direct proof of \eqref{YangBaxter}, based on the orthogonal 
decompositions of vector $X$
fields on $\g_{reg}^*$ resp. $G_{reg}$ into spherical and radial parts, 
$X=X^{sp}+ X^{rad}$. Here
the spherical part $X^{sp}$ is by definition tangent to orbits
and the radial part $X^{rad}$ orthogonal to orbits. Both radial and 
spherical vector fields are Lie-subalgebras of the Lie algebra of 
vector fields. 
It is convenient to introduce 
certain canonical skew-symmetric tensor fields
$\mathfrak{r}_0\in C^\infty(\g^*_{reg},\g\wedge\g)$ and
$\mathfrak{r}\in C^\infty(G_{reg},\g\wedge\g)$ on the set of regular
elements $\g^*_{reg}$ resp. $G_{reg}$.  Recall that $\g^*_{reg}$ is
the set of all $\mu\in\g^*_{reg}$ such that the stabilizer $G_\mu$
under the coadjoint action has minimal dimension (so that it is a
maximal torus of $G$). For all $\mu\in\g^*$ the operator $\ad_\mu$
(where we use the invariant inner product to identify $\g^*\cong\g$)
is invertible on the subspace $\g_\mu^\perp$. Let
$\ad_\mu^{-1}:\g\to\g$ be its extension to an operator on $\g$,
defined to be $0$ on $\g_\mu$. We set 
\begin{equation}\label{Defr0} 
\mathfrak{r}_0(\mu)_{ab}=\ad_\mu^{-1}(e_a)\cdot e_b.
\end{equation}
We remark that the definition of $\mathfrak{r}_0$ is in fact independent
of the inner product on $\g$. Viewing $\mathfrak{r}_0$ as a 2-form 
on $\g^*_{reg}$, it restricts to the Kirillov-Kostant-Souriau form 
on coadjoint orbits. 
The subset $G_{reg}$ consist of all $g\in G$ such that the stabilizer
$G_g$ is a maximal torus. Writing $g=\exp\mu$ the exponential map 
gives an isomorphism $\g_g^\perp=T_g(G\cdot g)\cong T_\mu(G\cdot\mu)\cong
\g_\mu^\perp\subset\g$. The skew-symmetric operator 
$\tanh(\f{\ad_\mu}{2})$ is invertible on $\g_\mu^\perp=\g_g^\perp$.
Let $\Big(\tanh(\f{\ad_\mu}{2})\Big)^{-1}$ be the extension of 
this inverse by $0$ to an operator on all of $\g$, and put
\begin{equation}\label{Defr}
\mathfrak{r}(g)_{ab}=\hh 
\Big(\tanh(\f{\ad_\mu}{2})\Big)^{-1}e_a\cdot e_b.
\end{equation}
Then 
\begin{equation}
T_{ab}=(\mathfrak{r}_0)_{ab}-\exp^*\mathfrak{r}_{ab}
\end{equation}
over $\exp^{-1}(G_{reg})\subset\g_{reg}$. 
Next, we examine the properties of the tensor fields
$\mathfrak{r}_0$ and $\mathfrak{r}$.
\begin{lemma}\label{r0}
The tensor field $\rr_0\in C^\infty(\g_{reg},\wedge^2\g)$ 
satisfies the equation 
$$ \on{Cycl}_{abc}\Big(\big(\f{\p}{\p\mu_a}\big)^{rad}(\rr_0)_{bc}
+(\rr_0)_{ak}f_{kbl}(\rr_0)_{lc}
\Big)=0.$$
\end{lemma}
\begin{proof}
We begin begin by observing that at any point $\mu\in\g_{reg}$, the
tensor $\on{Cycl}_{abc}(\ldots)$ on the left hand side of this
equation takes values in $\wedge^3\g_\mu^\perp$. This follows because
$\rr_0(\mu)$ takes values in $\wedge^2\g_\mu^\perp$, and the same is
true for radial derivatives
$\big(\f{\p}{\p\mu_a}\big)^{rad}(\rr_0)(\mu)$.   

We calculate
the spherical part of $ [\f{\p}{\p\mu_a},\f{\p}{\p\mu_b}]$ in two
ways. First, since partial derivatives commute we have 
$([\f{\p}{\p\mu_a},\f{\p}{\p\mu_b}])^{sp}=0$. On the other hand, 
decomposing into radial and spherical parts we have
$$
\big([\f{\p}{\p\mu_a},\f{\p}{\p\mu_b}]\big)^{sp}
=[(\f{\p}{\p\mu_a})^{sp},(\f{\p}{\p\mu_b})^{sp}]
+[(\f{\p}{\p\mu_a})^{rad},(\f{\p}{\p\mu_b})^{sp}]
+[(\f{\p}{\p\mu_a})^{sp},(\f{\p}{\p\mu_b})^{rad}].
$$
By definition of $\rr_0$, the spherical part of 
$\f{\p}{\p\mu_a}$ is $(\rr_0)_{ab}(e_b)_\g$. Since
$\rr_0$ is equivariant, 
$$ (e_a)_\g (\rr_0)_{bc}+f_{abk}(\rr_0)_{kc}+f_{acl}(\rr_0)_{bl}=0.$$
Hence the first term becomes 
\beq
[(\f{\p}{\p\mu_a})^{sp},(\f{\p}{\p\mu_b})^{sp}]
&=&[(\rr_0)_{ar}(e_r)_\g, \ (\rr_0)_{bs}(e_s)_\g]\\
&=&\on{Cycl}_{abc}\Big((\rr_0)_{ar}f_{rbs}(\rr_0)_{sc}\Big) (e_c)_\g.
\eeq
The second term and third term add up to 
$$
[(\f{\p}{\p\mu_a})^{rad},(\f{\p}{\p\mu_b})^{sp}]
+[(\f{\p}{\p\mu_a})^{sp},(\f{\p}{\p\mu_b})^{rad}]=
\Big((\f{\p}{\p\mu_a})^{rad}(\rr_0)_{bc}
+ (\f{\p}{\p\mu_b})^{rad}(\rr_0)_{ca}\Big)(e_c)_\g.
$$
Since $\big((\f{\p}{\p\mu_c})^{rad}(\rr_0)_{ab}\big)(e_c)_\g=0$
by orthogonality of radial and spherical vector fields, this can 
also be written
$$ \on{Cycl}_{abc}\big((\f{\p}{\p\mu_a})^{rad}(\rr_0)_{bc}\big)(e_c)_\g.$$
We have thus shown 
$$ \on{Cycl}_{abc}\Big(\big(\f{\p}{\p\mu_a}\big)^{rad}(\rr_0)_{bc}
+(\rr_0)_{ak}f_{kbl}(\rr_0)_{lc}
\Big)(e_c)_\g=0.$$
which proves the Lemma since the vector fields $(e_c)_\g$ span $\g_\mu^\perp$ 
for all $\mu\in\g_{reg}$.
\end{proof}

\begin{lemma}\label{r}
The tensor field $\rr \in C^\infty(G_{reg},\wedge^2\g)$ 
satisfies the equation 
$$ \on{Cycl}_{abc}\Big(\big(\f{e_a^L+e_a^R}{2}\big)^{rad}\rr_{bc}
+\rr_{ak}f_{kbl}\rr_{lc}
\Big)=\f{1}{4}f_{abc}.$$
\end{lemma}
\begin{proof}
Since $T_{ab}=(\rr_0)_{ab}-\exp^*\rr_{ab}$, Equation \eqref{Halfsum}
shows that  the spherical part of $\f{e_a^L+e_a^R}{2}$ is given by 
$$ \Big(\f{e_a^L+e_a^R}{2}\Big)^{sp}=\rr_{ab}(e_b)_G.$$
Using this Equation the proof of Lemma \ref{r} becomes parallel
to that of Lemma \ref{r0} --  the only difference being that 
$$ [\f{e_a^L+e_a^R}{2},\f{e_b^L+e_b^R}{2}]=\f{1}{4}f_{abc}(e_c)_G.$$
which accounts for the term $\f{1}{4}f_{abc}$ on the right hand side. 
\end{proof}
\begin{proof}[Proof of Equation \eqref{YangBaxter}]
As a consequence of
Equation \eqref{Halfsum},
the pull-back to $\Phinv(G_{reg})$ of the radial part of  
$\f{e_a^L+e_a^R}{2}$ is equal to the radial part of $\f{\p}{\p\mu_a}$.
Hence, combining Lemmas \ref{r0} and \ref{r} we find that 
\begin{equation}\label{radial}
 \on{Cycl}_{abc}\Big(\big(\f{\p}{\p\mu_a}\big)^{rad}T_{bc}
+ (\rr_0)_{ak}f_{kbl}(\rr_0)_{lc}
-\exp^* \rr_{ak}f_{kbl}\exp^*\rr_{lc}
\Big)=\f{1}{4}f_{abc}
\end{equation}
The equivariance property of $T$ shows that
\beq 
\big(\f{\p}{\p\mu_a}\big)^{sp}T_{bc}&=&
(\rr_0)_{ar}(e_r)_\g T_{bc}\\&=&
-(\rr_0)_{ar}(f_{rbs} T_{sc}+f_{rcs} T_{bs})\\
&=&
-2(\rr_0)_{ar} f_{rbs} (\rr_0)_{sc}+(\rr_0)_{ar} f_{rbs}\exp^* \rr_{sc}
+\exp^*\rr_{br} f_{rcs} (\rr_0)_{sa}.
\eeq
Taking the sum over cyclic permutations of $a,b,c$, adding to 
\eqref{radial}, and using $T_{ab}=(\rr_0)_{ab}-\exp^*\rr_{ab}$ 
 we obtain the identity \eqref{YangBaxter}.
\end{proof}

\end{appendix}

\bibliographystyle{plain}


\end{document}